\documentclass[11pt]{article}
\usepackage[utf8]{inputenc}

\usepackage{graphicx}
\graphicspath{ {./images/} }

\usepackage{standalone}

\usepackage{mathUtil}

\usepackage{fullpage}

\usepackage{comment}
\usepackage{indentfirst}

\usepackage{xcolor}

\usepackage[backend=biber,style=numeric,maxnames=50]{biblatex}
\usepackage{hyperref}
\hypersetup{
    colorlinks=true,
    linkcolor=blue,
    filecolor=magenta,
    urlcolor=cyan,
    pdftitle={Pi2-Rules and Inductive Classes of Godel Algebras},
    pdfpagemode=FullScreen,
}
\title{$\Pi_{2}$-Rule Systems and Inductive Classes of G\"{o}del Algebras}
\author{Rodrigo Nicolau Almeida}

\begin{document}

\maketitle

\begin{abstract}
    In this paper we present a general theory of $\Pi_{2}$-rules for systems of intuitionistic and modal logic. We introduce the notions of $\Pi_{2}$-rule system and of an inductive class, and provide model-theoretic and algebraic completeness theorems, which serve as our basic tools. As an illustration of the general theory, we analyse the structure of inductive classes of G\"{o}del algebras, from a structure theoretic and logical point of view. We show that unlike other well-studied settings (such as logics, or single-conclusion rule systems), there are continuum many $\Pi_{2}$-rule systems extending $\mathsf{LC}=\mathsf{IPC}+(p\rightarrow q)\vee (q\rightarrow p)$, and show how our methods allow easy proofs of the admissibility of the well-known Takeuti-Titani rule. Our final results concern general questions admissibility in $\mathsf{LC}$: (1) we present a full classification of those inductive classes which are inductively complete, i.e., where all $\Pi_{2}$-rules which are admissible are derivable, and (2) show that the problem of admissibility of $\Pi_{2}$-rules over $\mathsf{LC}$ is decidable.
\end{abstract}

\section{Introduction}

In the analysis of modal and intuitionistic logical systems, most research has focused on sets of theorems (in some literature called simply ``logics") and single-conclusion rule systems. Given a fixed language $\fancyL$, $\Gamma=\{\phi_{i} :i\leq n\}$ a finite (possibly empty) set of formulas and $\psi$ a formula in this language, we write
\begin{equation*}
\Gamma\vdash\psi
\end{equation*}
for a single-conclusion rule. In this setting, with such a rule we can associate a first-order formula
\begin{equation*}
\chi(\Gamma,\psi)\coloneqq \forall \overline{x} (\bigwedge_{i\leq n}\phi_{i}(\overline{x})\approx 1\rightarrow \psi(\overline{x})\approx 1)
\end{equation*}
such that an algebraic model $\alg{H}$ validates the rule $\Gamma\vdash\psi$ if and only if $\alg{H}\vDash \chi(\Gamma,\psi)$. These systems have been studied in a variety of settings, and a rich theory, connecting universal algebra and algebraic logic, has emerged out of their analysis (see e.g. \cite{Font2016-dk}), which covers several classical concepts, such as admissibility and structural completeness \cite{Iemhoff2015,DZIKstronkowskipaper,Cabrer2015}.

Nevertheless, more complex rules have also been considered for a long time: the pioneering work of Gabbay \cite{Gabbay1981} and Burgess \cite{Burgess1980}, in axiomatising temporal logics, and the work of Takeuti and Titani \cite{Takeuti1984-TAKIFL}, in axiomatisations of first-order G\"{o}del logic \cite{Baaz2017,metcalfemontagnasubstructuralfuzzy}, appear as early examples. Equally relevant is the work on ``Anti-Axioms" \cite{venemaantiaxioms}, and recent work axiomatising certain logics of compact Hausdorff spaces and studying structural properties of these calculi 
\cite{BEZHANISHVILI-stricimplicationcalculus,Bezhanishvili2022-if}. In all of these cases, $\Pi_{2}$-rules have been seen as an expedient way of obtaining axiomatisations for logical systems. Notwithstanding, such rules have not been the object of systematic study, and their model-theoretic and algebraic meaning has not been explored. This paper presents such an investigation, focusing on one hand on the development of a coherent theory of $\Pi_{2}$-rule systems and inductive rule classes, and on the other hand the admissibility of $\Pi_{2}$-rules for logical systems.

For the rest of the introduction we will provide an outline of the article, detailing the main ideas and contributions.  Throughout the whole paper we work with a given language $\fancyL$, which we will usually assume to be a language of intuitionistic or modal logic; we stress that such an assumption is not necessary, but makes the notation simpler, and will be sufficient for the examples at hand. We leave a general treatment of $\Pi_{2}$-rules, in full abstract algebraic generality, to further work.

\vspace{3mm}
\noindent \textbf{Definition 2.1} Let $\fancyL$ be a language. A $\Pi_{2}$-rule is a triple $(\Gamma,\psi,F)$ where:
\begin{enumerate}
    \item $\Gamma=\{\phi_{i}(\overline{p},\overline{q}): i\leq n\}$ and $\psi(\overline{p})$ are formulas in language $\fancyL$, the former possibly containing the variables in $\overline{p},\overline{q}$, and the latter containing at most variables in $\overline{p}$.
    \item $F=\{\overline{q}\}$ is a collection of propositional letters not occurring in $\psi$.
\end{enumerate}
Given such a triple, we often denote it by $\tilde{\forall}\overline{q} \ \Gamma(\overline{q})\vdash \psi$ or $\forall\overline{q}\Gamma(\overline{q})/^{2}\psi$, for short (sometimes without the universal quantifier, when the variables are clear from context or do not matter)\footnote{This notation serves to emphasise the fact both that these are distinct from usual rules, and the second-order nature of these rules, but we emphasise that it is purely formal. In particular, these rules do not involve propositional quantifiers.}:
\begin{prooftree}
\AxiomC{$\tilde{\forall} \overline{q}( \phi_{0}(\overline{p},\overline{q})\wedge...\wedge \phi_{n}(\overline{p},\overline{q}))$}
\UnaryInfC{$\psi(\overline{p})$}
\end{prooftree}
Given such a collection of formulas $\Gamma$, we refer to $F$ as the \textit{bound context of $\Gamma$}, and sometimes denote it as $B_{c}(\Gamma)$ when no confusion arises. If $F=\emptyset$ we write simply $\Gamma\vdash \psi$ and call this a $\Pi_{1}$-\textit{rule}. If $\Gamma=\emptyset$, we call this an \textit{axiom}.
\vspace{3mm}

In Section \ref{Pi2 rules and inductive rule classes} we introduce the notion of a $\Pi_{2}$-rule system, which provides a collection of proof rules for handling derivations with these rules. Given a fixed logical system $\mathsf{L}$, we can then form an extended calculus $\mathsf{L}\oplus \Sigma$ containing $\Pi_{2}$-rules. We then introduce an algebraic semantics for such rules, and prove the following completeness result, using a standard Lindenbaum-Tarski argument:

\vspace{3mm}
\noindent \textbf{Theorem 4.3} (Completeness Theorem for Inductive rules) Let $\Sigma$ be a collection of rules in the language $\fancyL$ Assume that $\Gamma/^{2}\phi\notin \mathsf{L}\oplus\Sigma$. Then there is some $\fancyL$-algebra $\alg{A}$, such that $\alg{A}\vDash \mathsf{L}\oplus\Sigma$, and $\alg{A}\nvDash \Gamma/^{2}\phi$.
\vspace{3mm}

In Section \ref{Model Theory of Inductive Rule classes}, we develop a theory paralleling that of varieties and quasivarieties for $\Pi_{2}$-rules. For that purpose, we recall that given a first-order formula $\phi$ (with terms in $\fancyL$), we say that $\phi$ is a $\forall\exists$-\textit{Special Horn formula} if it is of the form:
\begin{equation*}
    \forall\overline{x}(\forall\overline{y}(\bigwedge_{i=1}^{n}\lambda_{i}(\overline{x},\overline{y}))\rightarrow \gamma(\overline{x})),
\end{equation*}
where $\lambda_{i}$ and $\gamma$ are atomic formulas \footnote{Note that this is equivalent to the formula $\forall\overline{x}\exists\overline{y}(\bigwedge_{i=1}^{n}\lambda_{i}(\overline{x},\overline{y})\rightarrow \gamma(\overline{x}))$, which is where the terminology `$\Pi_{2}$' originates.}. Our objects of the study are thus the following classes:

\vspace{3mm}
\noindent \textbf{Definition 3.3} Let $\fancyL$ be a language, and $\alg{A},\alg{B}$ be two $\fancyL$-algebras. We say that $A$ is a $\forall$-\textit{subalgebra} of $B$, and write $\alg{A}\leq_{\forall}\alg{B}$ if it is a subalgebra, and for each atomic formula $\phi(\overline{x},\overline{y})$, and sequence of elements $\overline{a}\in A$:
    \begin{equation*}
        \alg{A}\vDash \forall \overline{x} \phi(\overline{x},\overline{a}) \implies \alg{B}\vDash \forall \overline{x}\phi(\overline{x},\overline{a}).
    \end{equation*}
    We define the notion of  $\forall$-embedding in the obvious way. Given a class of algebras $\mathbf{K}$, we write $\mathbb{S}_{\forall}(\mathbf{K})$ for the collection of $\forall$-subalgebras of elements of $\mathbf{K}$.
\vspace{3mm}

\vspace{3mm}
\noindent\textbf{Definition 3.6} Let $\mathbf{K}$ be a class of algebras. We say that $\mathbf{K}$ is an \textit{inductive rule class} if it is closed under $\forall$-subalgebras, products and ultraproducts.
\vspace{3mm}

We proceed to show in Theorem \ref{Collective theorem for special horn classes} that inductive rule classes are exaclty those axiomatised by Special Horn classes, through an adaptation of Mal'tsev's theorem. This is used in Section \ref{Model Theoretic Completeness for Pi2-rules} to provide a result connecting $\Pi_{2}$-rule systems and inductive rule classes:

\vspace{3mm}
\noindent \textbf{Theorem 3.10} There is a dual isomorphism, $Ind$, between the lattice of $\Pi_{2}$-rule systems, and the lattice of inductive rule classes of $\fancyL$-algebras.
\vspace{3mm}

This result fully entitles us to study $\Pi_{2}$-rule systems using the methods of universal algebra. It is from this lens that the problem of \textit{admissibility} of $\Pi_{2}$-rules is discussed. Under such a perspective, having a given calculus  we say that a $\Pi_{2}$-rule $\Gamma/^{2} \psi$ is \textit{admissible} over a calculus $\vdash$ if whenever some substitution instance of the premise can be derived, so can the conclusion. We moreover show how one can isolate the admissibility problem down to $\Pi_{2}$-rules, by introducing the notion of \textit{hereditarily admissible} rules -- a property which, on $\Pi_{1}$-rules, coincides with derivability, but which comes apart in the case of $\Pi_{2}$-rules. Our main result in this respect, following some ideas from \cite{Bezhanishvili2022-if}, connects this problem to the existence of model companions and model completions:

\vspace{3mm}
\textbf{Theorem 5.11}
\noindent
Let $\bf{K}$ be an inductive rule class, and suppose that it has a model companion $T^{*}$. Then a $\Pi_{2}$-rule $\Gamma/^{2}\psi$ is hereditarily admissible in $\bf{K}$ if and only if $T^{*}\vDash \Gamma/^{2}\psi$.
\vspace{3mm}

 The analysis of admissibility is here restricted to the case of $\mathsf{LC}$, the linear calculus of G\"{o}del and Dummett\footnote{Throughout the paper we refer to this often as G\"{o}del Logic, as is common practice in the literature.}, which provides a natural test-case: it is a well-studied logical system, possessing a very simple structure in the basic equational setting, e.g., the algebraic models of such a calculus, G\"{o}del algebras, are locally finite, and the lattice of extensions of $\mathsf{LC}$ is isomorphic to $(\omega + 1)^{*}$. In this context we provide some applications of our machinery:
 \begin{enumerate}
     \item In Corollary \ref{Density rule is admissible} we show that the density rule of Takeuti and Titani is admissible;
     \item In Theorem \ref{Characterization of the inductively complete rule classes} we provide a characterization of the inductive rule classes for which, in their dual $\Pi_{2}$-rule system, the notions of admissibility and derivability coincide;
     \item In Corollary \ref{Admissibility of Pi2-rules is effectively recognizable} we show that the admissibility problem for $\Pi_{2}$-rules is effectively recognizable over $\mathsf{LC}$\footnote{We note that in the time since this paper was submitted, the author has developed, together with Silvio Ghilardi, a method for recognising admissibility of $\Pi_{2}$-rules in a larger class of systems, which includes $\mathsf{LC}$ \cite{almeidaghilardiunificationsimplevariable}. Nevertheless, the proof is very different from the one included here, using instead a form of unification theory.}.
 \end{enumerate}
 
All of these applications follow from the universal algebraic methods previously developed, together with some auxiliary facts about the model theory of G\"{o}del algebras.

We conclude in Section \ref{Section: Conclusion} with some final remarks and future directions for this line of research.

\section{$\Pi_{2}$-rules and $\Pi_{2}$-Rule Systems}\label{Pi2 rules and inductive rule classes}

In this section we introduce the general notion of $\Pi_{2}$-rules, outline the algebraic semantic meaning of validating such a rule, and provide a proof system for these rules. Throughout we denote algebras by the variables $\mathcal{A},\mathcal{B},...$; to refer to the domains of these algebras we will use the corresponding letter in italics, $A,B,C,...$. Classes of algebras will be denoted by variables $\mathbf{K},\mathbf{T},...$. Propositional variables will be denoted by symbols $p,q,r$, and formula variables by symbols $\phi,\psi,\chi$; collections of formulas will usually be denoted by $\Gamma,\Delta,\Theta,....$. Variables for logics will usually be of the form $L,S,T$.

\begin{definition}
Let $\fancyL$ be a language. A $\Pi_{2}$-rule is a triple $(\Gamma,\psi,F)$ where:
\begin{enumerate}
    \item $\Gamma=\{\phi_{i}(\overline{p},\overline{q}): i\leq n\}$ and $\psi(\overline{p})$ are formulas in language $\fancyL$, the former possibly containing the variables in $\overline{p},\overline{q}$, and the latter containing at most variables in $\overline{p}$.
    \item $F=\{\overline{q}\}$ is a collection of propositional letters not occurring in $\psi$.
\end{enumerate}
Given such a triple, we often denote it by $\tilde{\forall}\overline{q} \ \Gamma(\overline{q})\vdash \psi$, sometimes as $\forall\overline{q}\Gamma(\overline{q})/^{2}\psi$, for short (sometimes without the universal quantifier, when the variables are clear from context or do not matter)\footnote{This notation serves to emphasise the fact both that these are distinct from usual rules, and the second-order nature of these rules, but we emphasise that it is purely formal. In particular, these rules do not involve propositional quantifiers.}, and sometimes presented in rule form as follows:
\begin{prooftree}
\AxiomC{$\tilde{\forall} \overline{q}( \phi_{0}(\overline{p},\overline{q})\wedge...\wedge \phi_{n}(\overline{p},\overline{q}))$}
\UnaryInfC{$\psi(\overline{p}).$}
\end{prooftree}
Given such a collection of formulas $\Gamma$, we refer to $F$ as the \textit{bound context of $\Gamma$}, and sometimes denote it as $B_{c}(\Gamma)$ when no confusion arises. If $F=\emptyset$ we write simply $\Gamma\vdash \psi$ and call this a $\Pi_{1}$-\textit{rule}. If $\Gamma=\emptyset$, we call this an \textit{axiom}.
\end{definition}

\begin{example}
    Let $\fancyL=(\wedge,\vee,\neg,0,1,\dia)$ be the language of modal algebras. Recall from \cite{Gabbay1981}, Gabbay's, ``Irreflexivity rule":
    \begin{prooftree}
        \AxiomC{$\tilde{\forall} p \ (\neg(p\rightarrow \dia p)\vee \phi)$}
        \UnaryInfC{$\phi,$}
    \end{prooftree}
    where $p$ does not occur in $\phi$. This was used to axiomatise irreflexive frames over temporal logics.
\end{example}

\begin{example}\label{Pi2rulesNIckandSilvioversion}
    More generally, in \cite{Bezhanishvili2022-if}, rules of the following form were considered:
    \begin{prooftree}
        \AxiomC{$\tilde{\forall}\overline{p} \ (F(\overline{p},\overline{q})\rightarrow \chi)$}
        \UnaryInfC{$G(\overline{q})\rightarrow\chi,$}
\end{prooftree}
where $\overline{p}$ does not occur in $\chi$ or $G$.
\end{example}

The following example has also been extensively studied, and was noted in \cite{Bezhanishvili2022-if} as not fitting in the existing frameworks:

\begin{example}\label{Takeuti-Titani Rule}
    Let $\fancyL=(\wedge,\vee,\rightarrow,0,1)$ be the language of intuitionistic logic. Recall the following rule by Takeuti and Titani, used for axiomatisation purposes in \cite{metcalfemontagnasubstructuralfuzzy} and \cite{Baaz2017}:
    \begin{prooftree}
        \AxiomC{$\tilde{\forall} r \ (g\rightarrow ((p\rightarrow r)\vee (r\rightarrow q)\vee c))$}
        \UnaryInfC{$g\rightarrow ((p\rightarrow q)\vee c).$}
    \end{prooftree}
    Such a rule is often called the ``Density rule", since it is valid on a linear Heyting algebra if and only if such an algebra is dense.
\end{example}

As the examples suggest, rather than being used as syntactic objects, $\Pi_{2}$-rules have often been used to capture specific semantic classes which cannot be defined using $\Pi_{1}$-rules. Hence we wish from the start to associate to these rules a particular semantics. As usual, we refer to an algebra $\alg{A}$ as an $\fancyL$-algebra if it is of type $\fancyL$. We denote by $\mathsf{Prop}$ a fixed but arbitrary set of propositional variables.

\begin{definition}
    Let $\alg{A}$ be an $\fancyL$-algebra and $\Gamma/^{2}\psi$ a $\Pi_{2}$-rule. Let $v:\mathsf{Prop}\to \alg{A}$ be a valuation over $\alg{A}$, and $(\alg{A},v)$ be a pair (called an \textit{algebraic model}). We say that $(\alg{A},v)$ makes this rule true, and write $$(\alg{A},v)\vDash \Gamma/^{2}\psi$$ if and only if: if for each $\forall \overline{q}\phi_{n}(\overline{p},\overline{q})\in \Gamma$, and each valuation $v'$ which is equal to $v$ up to variables in $B_{c}(\Gamma)$, $v'(\phi_{n}(\overline{p},\overline{q}))=1$; then $v(\psi(\overline{q}))=1$. 
    
    We write $\alg{A}\vDash \Gamma/^{2}\psi$ if for each model $(\alg{A},v)$ over $\alg{A}$, $(\alg{A},v)\vDash \Gamma/^{2}\psi$.    Given a collection of $\fancyL$-algebras $\bf{C}$, we write $\bf{C}\vDash \Gamma/^{2}\psi$ if $\alg{A}\vDash \Gamma/^{2}\psi$ for each $\alg{A}\in \bf{C}$.
\end{definition}

Note that in light of this definition, given any $\Pi_{2}$-rule there is a canonical way to assign to it a first-order formula. Namely, given $\Gamma(\overline{p},\overline{q})$ and $\psi(\overline{p})$, and $F=\{\overline{q}\}$ we let
\begin{equation*}
    \chi(\Gamma,\psi)=\forall \overline{y}(\forall \overline{x}(\bigwedge\Gamma(\overline{x},\overline{y})\approx 1) \rightarrow \psi(\overline{y})).
\end{equation*}
where $\Gamma=\{\phi_{0},...,\phi_{n}\}$ and $\bigwedge \Gamma(\overline{x},\overline{y})=\phi_{0}(\overline{x},\overline{y})\wedge...\wedge \phi_{n}(\overline{x},\overline{y})$. Then it is clear to see that:
\begin{equation*}
    \alg{H}\vDash \Gamma/^{2}\psi \iff \alg{H}\vDash \chi(\Gamma,\psi),
\end{equation*}
where the right hand side consequence is the usual first-order logic modelling relation. We will use this fact, often tacitly, throughout.

Having these rules, we can now come to the crucial definition of this section: a $\Pi_{2}$-Rule system. We recall some usual notation: given a formula $\chi$ containing variables $\overline{p}$, and $\overline{\phi}$ a sequence of formulas of the same length, we let $\chi[\overline{\phi}/\overline{p}]$ denote the uniform substitution of $\phi_{i}$ for $p_{i}$. We extend this notation to sets of formulas, i.e., $\Gamma[\overline{\phi}/\overline{p}]$, in the obvious way. Given a substitution $\sigma$, and $\Gamma$ a set of formulas, we write $\sigma[\Gamma]=\{\sigma(\psi) : \psi\in \Gamma\}$.

\begin{definition}\label{Definition of Pi2-rules}
    Let $\vdash$ be a set of $\Pi_{2}$-rules in the language $\fancyL$. We say that $\vdash$ is a $\Pi_{2}$-\textit{rule system} if it satisfies the following:
    \begin{itemize}
        \item (Strong Reflexivity) If $\overline{\phi}$ is a set of formulas not containing any variable in $\overline{p}$, we have $\tilde{\forall} \overline{p} \ \Gamma\vdash \chi[\overline{\phi}/\overline{p}]$.
        \item (Monotonicity) If $\tilde{\forall} \overline{p} \ \Gamma\vdash \psi$, then for any finite set of formulas $\Gamma'$, we have $\tilde{\forall} \overline{p} \ (\Gamma,\Gamma')\vdash \psi$.
        \item $(\alpha$-Renaming) if $\tilde{\forall}\overline{p} \ \Gamma\vdash \psi$, then $\tilde{\forall} \overline{q} \ \Gamma'\vdash \psi$, where $\Gamma'=\Gamma[\overline{q}/\overline{p}]$, so long as none of the variables from $\overline{q}$ occur free in either $\Gamma$ or $\psi$.
        \item (Rule cut) If  $\Gamma(\overline{p},\overline{r})$ and $\Delta(\overline{q},\overline{r})=\{\mu_{i}(\overline{q},\overline{r}) : i\leq n\}$, are two finite sets of formulas, $\psi$ does not contain any variables in $\overline{r}$ or $\overline{q}$, $\Gamma$ does not contain any variables from $\overline{q}$, and
        \begin{equation*}
            \tilde{\forall}\overline{p} \ \Gamma(\overline{p},\overline{r})\vdash \mu_{i}(\overline{q},\overline{r}) \text{ for each $i\leq n$ and } \tilde{\forall}\overline{q} \ \Delta(\overline{q},\overline{r})\vdash \psi,
        \end{equation*}
               
        Then we must have $\tilde{\forall}\overline{p} \ \Gamma(\overline{p},\overline{r})\vdash \psi$.
        \item (Bound Structurality) if $\tilde{\forall}\overline{p} \ \Gamma\vdash \psi$ and $\sigma$ is a substitution such that all variables in $F=\{\overline{p}\}$ are fixed, and such that $p\in F$ does not occur in $\sigma(q)$ for $q$ occuring free in $\Gamma$, $\sigma[\Gamma]$ or $\sigma(\phi)$, then $\tilde{\forall}\overline{p} \ \sigma[\Gamma]\vdash \sigma(\psi)$.
    \end{itemize}
\end{definition}

These definitions have some intuitive appeal deriving from our  notation. The following proposition establishes that the satisfaction relation for $\Pi_{2}$-rules forms a $\Pi_{2}$-rule system as above, and levies these ideas.

\begin{proposition}\label{Satisfaction of Pi-2 rules is a Pi-2 Heyting relation}
    Let $\alg{A}$ be an $\fancyL$-algebra. Then the set:
    \begin{equation*}
        \Pi_{2}(\alg{A})=\{\Gamma/^{2}\psi : \alg{A}\vDash \Gamma/^{2}\psi\}
    \end{equation*}
    forms a $\Pi_{2}$-rule system. More generally, given a class $\bf{K}$ of algebras, $\Pi_{2}(\bf{K})$ also forms a $\Pi_{2}$-rule system.
\end{proposition}
\begin{proof}
As noted above, we have that for each algebra $\alg{A}$, and each $\Pi_{2}$-rule $\Gamma/^{2}\psi$:
\begin{equation*}
    \alg{A}\vDash \Gamma/^{2}\psi \iff \alg{A}\vDash \chi(\Gamma,\psi).
\end{equation*}
With this in mind, one can check that the conditions of Definition \ref{Definition of Pi2-rules} follow from the rules for universal quantifiers of first-order logic. The result for classes of logics follows similarly.\end{proof}

With the former proposition in mind, we also note the following:

\begin{lemma}\label{Smallest rule system exists}
Let $\fancyL$ be a fixed language. If $S$ is a set of $\Pi_{2}$-rules, then there is a smallest $\Pi_{2}$-rule system containing $S$.
\end{lemma}
\begin{proof}
Note that the collection of $\Pi_{2}$-rules is closed under arbitrary intersections, since each of the conditions for being a $\Pi_{2}$-rule system is a closure condition; moreover, there is at least one such relation extending $S$ (namely, the relation containing all pairs of the form $\Gamma/^{2}\psi$).\end{proof}

\section{Inductive Rule Classes}\label{Model Theory of Inductive Rule classes}

Recall that given a first-order formula $\phi$, we say that $\phi$ is a $\forall\exists$-\textit{Special Horn formula} if it is of the form:
\begin{equation*}
    \forall\overline{x}(\forall\overline{y}(\bigwedge_{i=1}^{n}\lambda_{i}(\overline{x},\overline{y}))\rightarrow \gamma(\overline{x})),
\end{equation*}
where $\lambda_{i}$ and $\gamma$ are atomic first-order formulas. Such formulas, and some generalizations therein, were first studied by Lyndon \cite{lyndon1959properties}, in the context of seeking a characterisation of those formulas preserved under Horn formulas. They have moreover appeared several times in algebraic contexts, as well as in the algebra of logic; we give here two examples:

\begin{example} (Subordination Algebras)
Consider the language of Boolean algebras together with a binary symbol $\rightsquigarrow$. This is the language of so-called \textit{subordination algebras}. Then consider the following formula:
\begin{equation*}
    \forall a, b, d (\forall c (a\rightsquigarrow c \wedge c\rightsquigarrow b \leq d) \rightarrow a\rightsquigarrow b \leq d)
\end{equation*}
This formula was introduced in \cite{BEZHANISHVILI-stricimplicationcalculus}; it defines an important class of \textit{transitive subordination algebras}.\qed \end{example}

\begin{example} (Orthoimplicative systems)
Consider the language of orthoimplicative systems; this consists of the language of ortholattices (and Boolean algebras), together with infinitely many implications $(\multimap_{n})_{n\in \omega}$. Consider the following formula:
\begin{align*}
    \forall c_{0},...,c_{n},d_{0},...,d_{k},e &( (c_{i}=\bigvee c_{i}\wedge d_{j} \ \& \ \forall f ( f\wedge (\bigvee c_{i}\wedge d_{j}) \leq \bigvee (f\wedge c_{i}\wedge d_{j})))\\
    &\implies e\multimap_{n}(c_{0},...,c_{n})\leq e\multimap_{n}(d_{0},...,d_{k}) ) 
\end{align*}
This was shown in \cite{minhatesedemestrado} to be an interesting property of ortholattices, allowing a nicer relationship between ortholattices and modal algebras. Intuitively, it says that when the sequences of elements $c_{0},..,c_{n}$ and $d_{0},...,d_{k}$ form what is called an \textit{admissible join}, and the joins of the first sequence is below the second sequence, then the implication should respect this structure.\qed
\end{example}

From the discussion of Section \ref{Pi2 rules and inductive rule classes}, we can see that in model-theoretic terms, $\Pi_{2}$-rules amount to $\forall\exists$-Special Horn formulas. As such, in this section we will provide a number of results, analogous to Birkhoff's theorem on equational classes and Mal'tsev's theorem on quasi-equational classes. For that purpose, below, given a class of algebras $\mathbf{K}$, we will denote by:
\begin{equation}
    Th_{\forall\exists}^{H}(\mathbf{K})=\{\phi : \phi \text{ is $\forall\exists$ Special Horn}, \mathbf{K}\vDash \phi\},
\end{equation}
the inductive Special Horn consequences of $\mathbf{K}$. The following construction will play a key role in the whole theory:

\begin{definition}
    Let $\fancyL$ be a language, and $\alg{A},\alg{B}$ be two $\fancyL$-algebras. We say that $A$ is a $\forall$-\textit{subalgebra} of $B$, and write $\alg{A}\leq_{\forall}\alg{B}$ if it is a subalgebra, and for each atomic formula $\phi(\overline{x},\overline{y})$, and sequence of elements $\overline{a}\in A$:
    \begin{equation*}
        \alg{A}\vDash \forall \overline{x} \phi(\overline{x},\overline{a}) \implies \alg{B}\vDash \forall \overline{x}\phi(\overline{x},\overline{a}).
    \end{equation*}
    We define the notion of  $\forall$-embedding in the obvious way. Given a class of algebras $\mathbf{K}$, we write $\mathbb{S}_{\forall}(\mathbf{K})$ for the collection of $\forall$-subalgebras of elements of $\mathbf{K}$.
\end{definition}

Note that the latter operator\footnote{Technically, it is the operator $\mathbb{I}\mathbb{S}_{\forall}$ which is idempotent, but this will not affect anything.} is monotone, extensive and idempotent. We also recall the class operators $\mathbb{P},\mathbb{H}, \mathbb{S},\mathbb{P}_{S},\mathbb{P}_{U},\mathbb{I}$, respectively, closure under products, homomorphic images, subalgebras, subdirect products, ultraproducts, and isomorphisms, which we assume the reader to be familiar with. We use $\mathbb{V}$ to denote the closure under being a variety, i.e., $\mathbb{V}=\mathbb{HSP}$; and we use $\mathbb{Q}$ to denote closure under being a quasivariety, i.e., $\mathbb{Q}=\mathbb{ISPP}_{U}$. Given our usage of products and ultraproducts, whenever $a,b\in \prod_{i\in I}\alg{A}_{i}$, we use the following notation
\begin{equation*}
    \llbracket a=b\rrbracket =\{i\in I : a(i)=b(i)\}.
\end{equation*}
More generally, for a first-order formula $\phi(x_{1},...,x_{n})$, and elements $a_{1},...,a_{n}\in \prod_{i\in I}\alg{A}_{i}$ we write
\begin{equation*}
    \llbracket \phi(a_{1},...,a_{n})\rrbracket =\{i\in I : \alg{A}_{i}\vDash \phi(a_{1}(i),...,a_{n}(i))\}.
\end{equation*}
Moreover, if $U$ is an ultrafilter, we denote elements of $\prod_{i\in I}\alg{A}_{i}/U$, the ultraproduct, by $a/U$.

\begin{lemma}\label{Preservation of special Horn sentences}
If $\phi$ is a $\forall\exists$-Special Horn sentence, then $\phi$ is preserved under $\forall$-subalgebras, direct products and ultraproducts.
\end{lemma}
\begin{proof}
Assume that $\alg{A}\leq_{\forall}\alg{B}$, and we have a $\forall\exists$-Special Horn formula $\phi=\forall\overline{x}(\forall\overline{y}(\lambda_{0}(\overline{x},\overline{y})\wedge...\wedge \lambda_{n}(\overline{x},\overline{y}))\rightarrow \gamma(\overline{x}))$ which holds in $\alg{B}$. Now let $\overline{a}$ be any sequence in $A$. Assume that $\alg{A}\vDash \forall\overline{y}(\overline{\lambda}(\overline{a},\overline{y})$. Then by the assumption of this being a $\forall$-subalgebra, the same formula holds in $\alg{B}$. Since we assumed that $B$ satisfied the formula $\phi$, then $\alg{B}\vDash \gamma(\overline{a})$; but since $\alg{A}$ is a subalgebra, then $\alg{A}\vDash \gamma(\overline{a})$, as desired.

The fact that this sentence is preserved under direct products follows from general results on preservation of Horn sentences by direct products, whilst preservation by ultraproducts is true of any first-order formula.
\end{proof}

An important example which witnesses the usefulness of the concept of $\forall$-subalgebras is the case of subdirect products:

\begin{proposition}\label{Subdirect products are pi2subalgebras}
Let $\mathbf{K}$ be an arbitrary class of $\fancyL$-algebras. Then:
\begin{equation*}
    \mathbb{P}_{S}(\mathbf{K})\subseteq \mathbb{S}_{\forall}\mathbb{P}(\mathbf{K})
\end{equation*}
\end{proposition}
\begin{proof}
Assume that $f:\alg{A}\to \prod_{i\in I}\alg{B}_{i}$ is a map witnessing the fact that $\alg{A}$ is a subdirect product of the $\alg{B}_{i}$. Then we claim that $f$ $A$ is a $\forall$-embedding. Indeed, assume that $\alg{A}\vDash \forall \overline{x}\phi(\overline{x},\overline{a})$, where $\phi$ is an atomic formula. Suppose that $\prod_{i\in I}\alg{B}_{i}\nvDash \forall\overline{x}\phi(x,f(\overline{a}))$; hence for some tuple $\overline{b}$ of elements in $\prod_{i\in I}\alg{B}_{i}$,  $\prod_{i\in I}\alg{B}_{i}\nvDash \phi(\overline{b},f(\overline{a}))$. Thus for some $\alg{B}_{i}$, $B_{i}\nvDash \phi(\overline{b}(i),f(\overline{a})(i))$; since the product is subdirect, $\alg{B}_{i}$ is a homomorphic image of $\alg{A}$, so:
\begin{equation*}
    \alg{B}_{i}\nvDash \phi(f(\overline{c})(i),f(\overline{a})(i))
\end{equation*}
for some $c_{k}\in \alg{A}$ such that $f(c_{k})(i)=b_{k}(i)$, for each element of the tuple $\overline{b}$. But this in turn implies that $\prod_{i\in I}\alg{B}_{i}\nvDash \phi(f(\overline{c}),f(\overline{a}))$. Since this is an atomic formula, and $f$ is an embedding, we then have that $\alg{A}\nvDash \phi(\overline{c},\overline{a})$, which contradicts the hypothesis on $\alg{A}$. So by reductio, we have that $\alg{B}\vDash \forall\overline{x}\phi(x,f(a_{0}),...,f(a_{n}))$.\end{proof}

\begin{definition}
    Let $\mathbf{K}$ be a class of algebras. We say that $\mathbf{K}$ is an \textit{inductive rule class} if it is closed under isomorphic copies of $\forall$-subalgebras, products and ultraproducts.
\end{definition}

Note that as usual, the arbitrary intersection of inductive rule classes is again an inductive rule class. Hence we can consider an operator $\mathbb{R}^{\mathbb{I}}$, which constructs the smallest inductive rule class containing the given collection of algebras. As in the case of Birkhoff and Matlsev's results -- and owing to the proofs of the latter -- we also have a characterization of this operator:

\begin{theorem}\label{Characterisation of inductive rule classes}
Let $\mathbf{K}$ be an arbitrary collection of algebras. Then we have that:
\begin{equation*}
    \mathbb{R}^{\mathbb{I}}(\mathbf{K})=\mathbb{IS}_{\forall}\mathbb{PP}_{U}(\mathbf{K})
\end{equation*}
\end{theorem}
\begin{proof}
Since inductive rule classes are closed under isomorphic copies of $\forall$-subalgebras, products and ultraproducts, we must have that $\mathbb{IS}_{\forall}\mathbb{PP}_{U}(\mathbf{K})\subseteq \mathbb{R}^{\mathbb{I}}(\mathbf{K})$. For the other inclusion, note first that:
\begin{equation*}
    \mathbb{P}\mathbb{S}_{\forall}\leq \mathbb{S}_{\forall}\mathbb{P}
\end{equation*}
Indeed if $\alg{A}=\prod_{i\in I}\alg{B}_{i}$, where each of the factors is a $\forall$-subalgebra of $\alg{C}_{i}$, then $\alg{A}$ will be a $\forall$-subalgebra of $\prod_{i\in I}\alg{C}_{i}$; it is surely a subalgebra, and if $\forall \overline{x}\phi(\overline{x},\overline{a})$ is true in $\alg{A}$, then $\forall \overline{x}\phi(\overline{x},\overline{a}(i))$ is true in $\alg{B}_{i}$ for each $i$, hence $\forall \overline{x}\phi(\overline{x},\overline{a}(i))$ is true in $\alg{C}_{i}$, so $\forall \overline{x}\phi(\overline{x},\overline{a})$ is true in $\prod_{i\in I}\alg{C}_{i}$.

Similarly, we have that \begin{equation*}
    \mathbb{P}_{U}\mathbb{S}_{\forall}\leq \mathbb{S}_{\forall}\mathbb{P}_{U}
\end{equation*}
Indeed, if $\alg{A}$ is an ultraproduct of $\alg{B}_{i}$ which are $\forall$-subalgebras
 of $\alg{C}_{i}$ via $g_{i}$, then consider $\prod_{i\in I}\alg{C}_{i}$, and take the ultraproduct by the same ultrafilter. Let $g:\prod_{i\in I}\alg{B}_{i}\to \prod_{i\in I}\alg{C}_{i}$ be the composite function. Define the map $f:\alg{A}\to \prod_{i\in I}C_{i}/U$ as $f(a/U)=g(a)/U$. Then this is well-defined and injective: if $a/U\neq b/U$, then $I-\llbracket a=b\rrbracket \in U$, so $I-\llbracket g(a)=g(b)\rrbracket \in U$ (because $g$ is injective, and $U$ is upwards closed), so $g(a)/U\neq g(b)/U$; the well-definedness condition follows by dual arguments. It is a homomorphism by usual arguments. Moreover, if $\phi(\overline{x},\overline{a})$ is any formula such that $\alg{A}\vDash \forall \overline{x}\phi(\overline{x},\overline{a})$, then $\llbracket \forall \overline{x}\phi(\overline{x},\overline{a})\rrbracket \in U$; for each $i\in \llbracket \forall \overline{x}\phi(\overline{x},\overline{a})\rrbracket $, then $\alg{B}_{i}\vDash \forall \overline{x}\phi(\overline{x},\overline{a}(i))$, so by assumption, $\alg{C}_{i}\vDash \forall \overline{x}\phi(\overline{x},g_{i}(\overline{a}(i)))$. Hence:
 \begin{equation*}
     \llbracket \forall \overline{x}\phi(\overline{x},g(\overline{a}))\rrbracket \in U
 \end{equation*}
Which shows the result, by \L o\'{s}'s theorem. 

Now as noted above, we have that $\mathbb{IS}_{\forall}(\mathbb{IS}_{\forall}\mathbb{PP}_{U}(\mathbf{K}))=\mathbb{IS}_{\forall}\mathbb{PP}_{U}(\mathbf{K})$. Additionally, note that:
\begin{equation*}
    \mathbb{P}(\mathbb{IS}_{\forall}\mathbb{PP}_{U}(\mathbf{K}))\subseteq \mathbb{IS}_{\forall}(\mathbb{PPP}_{U}(\mathbf{K}))=\mathbb{IS}_{\forall}\mathbb{PP}_{U}(\mathbf{K})
\end{equation*}
and also:
\begin{equation*}
    \mathbb{P}_{U}(\mathbb{IS}_{\forall}\mathbb{PP}_{U}(\mathbf{K}))\subseteq \mathbb{IS}_{\forall}\mathbb{P}_{U}\mathbb{PP}_{U}(\mathbf{K}))
\end{equation*}
But now notice that $\mathbb{P}_{U}\mathbb{P}\leq \mathbb{P}_{R}\mathbb{P}_{R}\leq \mathbb{P}_{R}$. Additionally, by \cite[Theorem 1.2.12]{Gorbunov1998-qd}, we have that $\mathbb{P}_{R}\leq \mathbb{P}_{S}\mathbb{P}_{U}$. Also by our Proposition \ref{Subdirect products are pi2subalgebras}, we have that $\mathbb{P}_{S}\leq \mathbb{S}_{\forall}\mathbb{P}$, and hence:
\begin{equation*}
    \mathbb{P}_{U}\mathbb{P}\leq \mathbb{IS}_{\forall}\mathbb{P}\mathbb{P}_{U}
\end{equation*}
This thus entails above that:
\begin{equation*}
    \mathbb{P}_{U}(\mathbb{IS}_{\forall}\mathbb{PP}_{U}(\mathbf{K}))\subseteq \mathbb{IS}_{\forall}\mathbb{PP}_{U}(\mathbf{K}),
\end{equation*}
as desired. Hence $\mathbb{IS}_{\forall}\mathbb{PP}_{U}(\mathbf{K})$ is closed under all operations, and is contained in the smallest inductive rule class, implying they must be equal.\end{proof}

\begin{definition}
    Let $\alg{A}$ be an $\fancyL$-algebra. We denote the $\forall$-\textit{open diagram} of $\alg{A}$ as follows:
    \begin{equation*}
        T_{0}(A)=\{\forall \overline{x}\phi(\overline{x},\overline{a}) : \phi(\overline{x},\overline{a}) \text{ is atomic }, \alg{A}\vDash \forall \overline{x}\phi(\overline{x},\overline{a}) \}\cup \{\chi(\overline{a}) : \chi(\overline{a}) \text{ is negated atomic}, \alg{A}\vDash \chi(\overline{a}) \},
    \end{equation*}
    where these formulas are considered in the language with constants from $A$.
\end{definition}

The following lemma will be very useful, and runs along the same lines as usual proofs using the method of diagrams:

\begin{lemma}\label{Diagrams for pi2-rules}
For each $\fancyL$-algebras $\alg{A},\alg{B}$ the following are equivalent:
\begin{enumerate}
    \item $\alg{A}\leq_{\forall}\alg{B}$;
    \item $\alg{B}$ is a model of $T_{0}(\alg{A})$.
\end{enumerate}
\end{lemma}

Putting all of this together gives us the main result of this section:

\begin{theorem}\label{Collective theorem for special horn classes}
Let $\mathbf{K}$ be a collection of algebras. Then the following are equivalent:
\begin{enumerate}
    \item $\mathbf{K}$ is an inductive rule class;
    \item $\mathbf{K}$ is axiomatised using $\forall\exists$ Special Horn formulas;
    \item $\mathbf{K}=\mathbb{IS}_{\forall}\mathbb{PP}_{U}(\mathbf{K}^{*})$ for some class $\mathbf{K}^{*}$.
\end{enumerate}
\end{theorem}
\begin{proof}
(3) is equivalent to (1) by Theorem \ref{Characterisation of inductive rule classes}. (2) implies (1) by Lemma \ref{Preservation of special Horn sentences}. So we show that (1) implies (2). For that purpose, we will consider $Th_{\forall\exists}^{H}(\bf{K})$, and show that this axiomatises $\bf{K}$. So assume that $\alg{A}$ is a $\fancyL$-algebra satisfying this set of formulas.

Now consider $\{\forall \overline{x}\phi_{1}(\overline{a},\overline{x}),...,\forall \overline{x}\phi_{n}(\overline{a},\overline{x}),\neg\chi_{0}(\overline{a}),...,\neg\chi_{k}(\overline{a})\}\subseteq T_{0}(A)$, where all formulas are atomic; then $$\mathcal{A}\vDash \exists \overline{a}(\forall \overline{x}(\phi_{1}(\overline{x},\overline{a})\wedge...\wedge \phi_{n}(\overline{x},\overline{a}))\wedge \neg\chi_{0}(\overline{a})\wedge...\wedge\neg\chi_{k}(\overline{a})).$$

Note that this follows from the fact that universal quantifiers commute with conjunction, and hence we can pull them outside of the conjunction of the $\phi_{i}$. We will show that some member of $\mathbb{P}(\mathbf{K})$ satisfies this sentence as well. For this purpose it suffices to show that:
\begin{equation*}
\mathbb{P}(\mathbf{K}) \nvDash \forall \overline{y} (\forall \overline{x}( \phi_{1}(\overline{y},\overline{x})\wedge...\wedge \phi_{n}(\overline{y},\overline{x})) \rightarrow \chi_{0}(\overline{y})\vee ... \vee \chi_{k}(\overline{y}))
\end{equation*}
If there were at most one $\chi_{i}$, then this would be logically equivalent to a $\forall\exists$ Special Horn formula which is not true of $A$, and since the latter satisfies all $\forall\exists$ Special Horn consequences of the theory of $\mathbf{K}$, not of $\mathbf{K}$. So let us suppose that there are at least two  $\chi_{i}$. Then for $1\leq i\leq k$ one can take the formula:
\begin{equation*}
    \mu_{i}\coloneqq \forall\overline{y}(\forall\overline{x}(\phi_{1}(\overline{y},\overline{x})\wedge...\wedge \phi_{n}(\overline{y},\overline{x}))\rightarrow \chi_{i}(\overline{y})).
\end{equation*}
And we will obtain that $\mathbf{K}\nvDash \mu_{i}$, since otherwise, we would have a contradiction to the fact that $\alg{A}$ satisfies all $\forall\exists$-Special Horn consequences of $\bf{K}$. Hence for some $\alg{A}_{i}\in \mathbf{K}$, $\alg{A}_{i}\vDash \exists \overline{y} \forall \overline{x} (\phi_{1}(\overline{x},\overline{y})\wedge \phi_{n}(\overline{x},\overline{y}))\wedge \neg\chi_{i}(\overline{y}))$. Now choose elements $a_{1}(i),...,a_{n}(i)\in \alg{A}_{i}$, such that:
\begin{equation*}
    \alg{A}_{i}\vDash \forall\overline{x}(\phi_{1}(\overline{x},a_{1}(i),...,a_{n}(i))\wedge ...\wedge \phi_{n}(\overline{x},a_{1}(i),...,a_{n}(i)))\wedge \neg\chi_{i}(a_{1}(i),...,a_{n}(i))
\end{equation*}
Then look at $\alg{A}_{0}\times...\times \alg{A}_{k}$, and consider the elements $a_{1},...,a_{k}$ defined to coincide with the witnesses above on each $i$. Then we have:
\begin{equation*}
    \alg{A}_{0}\times...\times \alg{A}_{k}\vDash \bigwedge_{i=1}^{n}\forall\overline{x}(\phi_{i}(\overline{x},a_{1},...,a_{k}) \wedge \bigwedge_{j=1}^{k}\neg\chi_{i}(a_{1},...,a_{k})
\end{equation*}
Which follows from the facts we had about $\forall\exists$-Special Horn formulas. Hence, we have shown the desired fact about $\mathbb{P}(\bf{K})$.

Now let $I$ be the set of finite subsets of $T_{0}(\alg{A})$. By the argument we just showed, for each $i\in I$, there is some $\alg{A}_{i}\in \mathbb{P}(\mathbf{K})$, and elements $\overline{a}(i)\in \alg{A}_{i}$ such that the formulas in $i$ become true of $\alg{A}_{i}$ when $\overline{a}$ is interpreted as $\overline{a}(i)$.

Hence consider $J_{i}=\{j\in I : i\subseteq j\}$. Let $F={\uparrow}\{J_{i} : i\in I\}$; then note that $F$ is a filter: it is clearly uwpards closed, and if $J_{i_{0}}\subseteq K$ and $J_{i_{1}}\subseteq K'$, note that $J_{i_{0}}\cap J_{i_{1}}=\{j\in I : i_{0}\subseteq j \text{ and } i_{1}\subseteq j\} = J_{i_{0}\cup i_{1}}$ so $J_{i_{0}\cup i_{1}}\subseteq K\cap K'$. Hence by the ultrafilter principle, let $U$ be an ultrafilter extending $F$.

Now let $\hat{a}$ be the tuple of elements in $\prod_{i\in I}A_{i}$ whose ith coordinate is $\overline{a}(i)$. Then  note that for $\psi(a_{0},...,a_{n})\in T_{0}(\mathcal{A})$, we have:
\begin{equation*}
\llbracket \psi(\hat{a}_{0},...,\hat{a}_{n})\rrbracket \supseteq J_{i}\in U
\end{equation*}
where $i=\{\psi(a_{0},...,a_{n})\}$; hence $\llbracket \psi(\hat{a}_{0},...,\hat{a}_{n})\rrbracket \in U$, so $\prod_{i\in I}\mathcal{A}_{i}/U\vDash \psi(\hat{a}_{0}/U,...,\hat{a}_{n}/U)$.

By hypothesis on the factors of the product, $\prod_{i\in I}\mathcal{A}_{i}/U\in \mathbb{P}_{U}\mathbb{P}(\mathbf{K})$. By what we just showed, this satisfies $T_{0}(\mathcal{A})$, so by Lemma \ref{Diagrams for pi2-rules}, we have that $\mathcal{A}\in \mathbb{I}\mathbb{S}_{\forall}\mathbb{P}_{U}\mathbb{P}(\mathbf{K})$. By the closure conditions, we then obtain that:
\begin{equation*}
    \mathcal{A}\in \mathbf{K},
\end{equation*}
which shows the axiomatisation, as intended.\end{proof}

We conclude this section by providing a practical criterion for being a $\forall$-embedding.

\begin{proposition}\label{Criterion for being a universal subalgebra}
Assume that $f:\alg{A}\to \alg{B}$ is an embedding. Then $f$ is a $\forall$-embedding if and only if for every finite sequence $\overline{a}$, $(\alg{B},f\overline{a})\in \mathbb{V}(\alg{A},\overline{a})$.
\end{proposition}
\begin{proof}
First assume that $f$ is not a $\forall$-embedding. Hence there is some $\overline{a}\in A$, and some equation $\phi(\overline{a},\overline{x})$ such that $A\vDash \forall \overline{x}\phi(\overline{a},\overline{x})$ and $B\nvDash \forall\overline{x}\phi(f\overline{a},\overline{x})$. But then clearly $(\alg{ B},f\overline{a})$ cannot belong to the variety generated by $(\alg{A},\overline{a})$: the equation
\begin{equation*} \phi(\overline{a},\overline{x})
\end{equation*}
is valid in $(\alg{A},\overline{a})$, but the same equation, interpreted over $(\alg{B},f\overline{a})$ fails for some $\overline{x}$; so by Birkhoff's theorem, we have that $(\alg{B},f\overline{a})\notin \mathbb{V}(\alg{A},\overline{a})$.

Conversely, assume that $\alg{A}\leq_{\forall}\alg{B}$. Then whenever $(\alg{A},\overline{a})$ validates some equation, so does $(\alg{B},f\overline{a})$ by assumption; hence by Birkhoff's theorem, $(\alg{B},f\overline{a})$ belongs to the variety generated by $(\alg{A},\overline{a})$.
\end{proof}

We will have opportunity to use this criterion often in Section \ref{Inductive rule classes of godel algebras}. For now, we can extract a simple corollary of this about Boolean algebras:

\begin{corollary}
    There exists a unique inductive rule class of Boolean algebras.
\end{corollary}
\begin{proof}
Let $\mathbf{B}$ be an arbitrary Boolean algebra. It is clear that $\mathbf{2}$ is a subalgebra of $\mathbf{B}$; we claim that moreover it is a universal subalgebra. Indeed, note that $\mathbf{2}$ admits no proper extension with constants (since all of its elements are already constants in the language). Additionally, $\mathbf{B}$ will certainly belong to the variety generated by $\mathbf{2}$. So by Proposition \ref{Criterion for being a universal subalgebra}, $\mathbf{2}$ is indeed a universal subalgebra. 

So notice that whenever $\mathbf{K}$ is a non-empty inductive rule class, $\mathbf{2}\in \mathbf{K}$; since the latter is closed by subdirect products by our results above, then $\mathbb{P}_{S}(\mathbf{2})\subseteq \mathbf{K}$; but by Stone's theorem, the former already contains all Boolean algebras.\end{proof}

\section{Algebraic Completeness for $\Pi_{2}$-rule Systems}\label{Model Theoretic Completeness for Pi2-rules}

The results of the previous section can in some sense be called a form of model-theoretic completeness. However, this is not the same as completeness with respect to the notion of $\Pi_{2}$-rule system we have previously developed. In this section, we show, using a Lindenbaum-Tarski style argument, that this can furthermore be made into a dual isomorphism between the lattice of $\Pi_{2}$-rule Systems and the lattice of inductive rule classes. Our results in this section follow the usual recipe for establishing such results, and the notion of proofs using $\Pi_{2}$-rules is generalised from \cite{BEZHANISHVILI-stricimplicationcalculus} and \cite{Bezhanishvili2022-if}.

Throughout this section we assume that we have a logical calculus $\mathsf{L}$ of normal modal logic, or intuitionistic logic, which is algebraized by some class of (Heyting or modal) algebras $\bf{K}$\footnote{The reader interested in more complex axiomatisations should not have difficulty modifying Definition 4.1 below: whatever logic $\mathsf{L}$ one takes, understood as a collection of $\Pi_{1}$-rules, one should simply add all the $\Pi_{1}$-rules to the possible derivations.}. If $\Sigma$ is a collection of $\Pi_{2}$-rules, we will denote by $\mathsf{L}\oplus \Sigma$ the smallest $\Pi_{2}$-rule system extending $\mathsf{L}$ which contains the rules from $\Sigma$ (which exists by Lemma \ref{Smallest rule system exists}).

\begin{definition}\label{Derivation using pi2-rules}
Let $\fancyL$ be an algebraic language, and let $\Sigma$ be a set of $\Pi_{2}$-rules. Let $\mathsf{L}$ be a logical calculus. If $(\Gamma,\phi,F)$ is a $\Pi_{2}$-rule, we say that $\phi$ is  \textit{derivable} from $\Gamma$ in $\mathsf{L}$ using the $\Pi_{2}$-rules in $\Sigma$, and write $\tilde{\forall}\overline{p} \ \Gamma\vdash_{\mathsf{L}\oplus \Sigma}\phi$, provided there is a sequence $\psi_{0},...,\psi_{n}$ of formulas such that:
\begin{itemize}
\item $\psi_{n}=\phi$;
\item For each $\psi_{i}$ we have that either:
\begin{enumerate}
\item $\psi_{i}$ is an instance of $\Gamma$ where possibly some free variables are replaced by some substitution instance;
\item $\psi_{i}$ is an instance of an axiom of $\mathsf{L}$ or,
\item $\psi_{i}$ is obtained using Modus Ponens from some previous $\psi_{j}$ or,
\item $\psi_{i}= \chi(\overline{\xi})$ and $\psi_{j_{k}}= \mu_{k}(\overline{r},\overline{\xi})$ for $j_{k}<i$, where
\begin{itemize}
    \item $\Delta=\{\mu_{k}(\overline{q},\overline{p}) : k\in \{0,...,m\}\}$;
    \item $\Delta/^{2}\chi\in \Sigma$;
    \item $\chi=\chi(\overline{p})$;
    \item $\overline{r}$ does not appear free in $\Gamma$.\footnote{Note that this definition works implicitly up to $\alpha$-renaming, which is done to avoid complicated tracking of substitutions renaming variables. All of this could be developed using such artifacts, though we will not pursue it here.}
\end{itemize}
\end{enumerate}
\end{itemize}
\end{definition}

\begin{lemma}
    Let $\Sigma$ be a collection of rules in the language $\fancyL$, and $\tilde{\forall}\overline{p} \ \Gamma /^{2}\phi$ one such rule. Then:
    \begin{equation*}
        \tilde{\forall}\overline{p} \ \Gamma\vdash_{\mathsf{L}\oplus \Sigma}\phi \iff \tilde{\forall}\overline{p} \ \Gamma/^{2}\phi\in \mathsf{L}\oplus\Sigma
    \end{equation*}
\end{lemma}
\begin{proof}
    Right to left is immediate. Left to right follows by an induction on the structure of the derivation: Strong Reflexivity ensures condition (1); bound structurality ensures that the substitutions work; from rule cut we get both the use of Modus Ponens, and allows us to pull through the case for other inductive rules being used in the derivation.
\end{proof}

\begin{theorem}[Completeness Theorem for Inductive rules]\label{Completeness Theorem for Inductive Rules}
    Let $\Sigma$ be a collection of rules in the language $\fancyL$ Assume that $\Gamma/^{2}\phi\notin \mathsf{L}\oplus\Sigma$. Then there is some $\fancyL$-algebra $\alg{A}$, such that $\alg{A}\vDash \mathsf{L}\oplus\Sigma$, and $\alg{A}\nvDash \Gamma/^{2}\phi$.
\end{theorem}
\begin{proof}
Since $\Gamma/^{2}\phi\notin \Sigma$, by the previous Lemma, $\tilde{\forall}\overline{p} \ \Gamma\nvdash_{\mathsf{L}\oplus \Sigma}\phi$.  Hence, let $\mathsf{Prop}$ be the set of proposition letters. Let $\freeP$ be the Lindenbaum-Tarski algebra defined over the term algebra, in the following way
\begin{equation*}
    [\psi]\leq [\mu]\iff  \tilde{\forall}\overline{p} \ \Gamma\vdash_{\mathsf{L}\oplus\Sigma} \psi\rightarrow\mu
\end{equation*}
The standard argument then shows that $\freeP$ yields a $\fancyL$-algebra. We additionally claim that $\freeP$ validates the axioms of $\Sigma$. Indeed take $\Delta/^{2}\psi$ a certain rule, and $v:\mathsf{Prop}\to \freeP$ a valuation, with the property that for each valuation $v'$ differing from it at most in $B_{c}(\Delta)$, the bound context of $\Delta$,  $v'(\Delta)=[\top]$. Note that then, in particular, we can take this for some propositional variables not ocurring free in $\Gamma$, since only finitely many do so; hence set $v'(q_{i})=q_{i}$ for $q_{i}$ not occuring free in $\Gamma$, and not ocurring in $v(p)$ for any $p$ which is not in $B_{c}(\Delta)$, and $v'(p)=v(p)$ otherwise, and consider $v'$ as a substitution, i.e. for each $p\notin B_{c}(\Delta)$, $v'(p)=\chi_{p}$ is some formula. Then we note that:
\begin{equation*}
    [\psi_{i}(\chi_{p},\overline{q})]=[\top]
\end{equation*}
It follows that for each $i$:
\begin{equation*}
    \tilde{\forall}\overline{p} \ \Gamma\vdash_{\mathsf{L}\oplus\Sigma} \psi_{i}(\overline{\chi_{p}},\overline{q})
\end{equation*}
Hence note that since $\Delta/^{2}\psi$ is in the calculus (and there are no variable clashes), we can then derive:
\begin{equation*}
    \tilde{\forall}\overline{p} \ \Gamma\vdash_{\mathsf{L}\oplus\Sigma} \psi(\overline{\chi_{p}}).
\end{equation*}
Which means that $v(\psi(\overline{p}))=1$. This shows that $\freeP\vDash \Delta/^{2}\psi$. We conclude that $\freeP\vDash \mathsf{L}\oplus \Sigma$.

Now consider the valuation over $\freeP$ defined by $v(p)=p$. Note that since $\Gamma\vdash\mu$ for each $\mu\in \Gamma$, we have that $[\bigwedge\Gamma]=1$. Additionally, for each valuation $v'$ differing at most outside $B_{c}(\Gamma)$, we can consider $v'$ as a substitution, and obtain that $\Gamma\vdash_{\mathsf{L}\oplus\Sigma}v'(\mu)$ again, by our assumption on the derivations. On the other hand, by hypothesis:
\begin{equation*}
    \Gamma\nvdash_{\mathsf{L}\oplus\Sigma}\phi
\end{equation*}
and so $v(\phi)\neq [\top]$. Hence, $\freeP,v\nvDash\Gamma/^{2}\phi$, as desired.\end{proof}

Having this completeness theorem, we now have all the tools necessary to prove the dual isomorphism as desired. The following definition was mentioned above:

\begin{definition}
    Let $\mathbf{K}$ be a class of $\fancyL$-algebras. Define the collection of $\Pi_{2}$-rules associated with $\mathbf{K}$ as follows:
    \begin{equation*}
        \Pi_{2}(\mathbf{K})\coloneqq \{ \Gamma/^{2}\psi : \alg{A}\vDash \Gamma/^{2}\psi, \forall \alg{A}\in \mathbf{K}\}
    \end{equation*}
\end{definition}

The following is a straightforward consequence of Proposition \ref{Satisfaction of Pi-2 rules is a Pi-2 Heyting relation}:

\begin{proposition}
    For each class $\mathbf{K}$ of $\fancyL$-algebras, $\Pi_{2}(\mathbf{K})$ is a $\Pi_{2}$-rule system.
\end{proposition}

As before we have that that the operations defining an inductive rule class preserve satisfaction of these rules, essentially given that such rules correspond to $\forall\exists$-Special Horn formulas:

\begin{proposition}
Let $\mathbf{K}$ be a class of $\fancyL$-algebras. Then:
\begin{equation*}
    \Pi_{2}(\mathbf{K})=\Pi_{2}(\mathbb{S}_{\forall}(\mathbf{K}))=\Pi_{2}(\mathbb{P}(\mathbf{K}))=\Pi_{2}(\mathbb{P}_{U}(\mathbf{K}))
\end{equation*}
Consequently:
\begin{equation*}
    \Pi_{2}(\mathbf{K})=\Pi_{2}(\mathbb{R}^{\mathbb{I}}(\mathbf{K}))
\end{equation*}
\end{proposition}

Likewise, we have the following operation, moving in the dual direction:

\begin{definition}
    Given a collection $\Sigma$ of $\fancyL$ $\Pi_{2}$-rules, define:
    \begin{equation*}
        Ind(\mathsf{L}\oplus \Sigma)\coloneqq \{ \alg{A} : \alg{A}\vDash \Gamma/^{2}\psi, \tilde{\forall}\overline{p} \  \Gamma/^{2}\psi\in \mathsf{L}\oplus \Sigma\}
    \end{equation*}
    We call this the \textit{inductive rule class} generated by $\mathsf{L}\oplus \Sigma$.
\end{definition}

We can now prove the following:

\begin{proposition}\label{Model Theoretic Part of completeness}
    For each inductive rule class of $\fancyL$-algebras $\mathbf{K}$, we have that:
    \begin{equation*}
        Ind(\Pi_{2}(\mathbf{K}))=\mathbf{K}
    \end{equation*}
\end{proposition}
\begin{proof}
One inclusion is clear: if $\alg{A}\in \mathbf{K}$, then it certainly validates all rules validated by all algebras in that class.
We focus on the other. By Theorem \ref{Collective theorem for special horn classes}, $\bf{K}$ is axiomatised by $\forall\exists$-formulas. So if $\alg{A}\in Ind(\Pi_{2}(\bf{K}))$, then $\alg{A}\vDash Th_{\forall\exists}^{H}(\bf{K})$, and so $\alg{A}\in \bf{K}$.\end{proof}

The other direction follows from our algebraic completeness proof:

\begin{corollary}\label{Proof Theoretic Part of Completeness}
    Let $\Sigma$ be a collection of $\fancyL$-rules. Then:
    \begin{equation*}
        \Pi_{2}(Ind(\mathsf{L}\oplus\Sigma))=\mathsf{L}\oplus\Sigma.
    \end{equation*}
\end{corollary}
\begin{proof}
    The right to left inclusion surely holds. Now assume that $\Gamma/^{2}\phi\notin \mathsf{L}\oplus \Sigma$. By Theorem \ref{Completeness Theorem for Inductive Rules}, we have that there is an algebra $\alg{A}$ such that $\alg{A}\vDash \mathsf{L}\oplus \Sigma$ and $\alg{A}\nvDash \Gamma/^{2}\phi$. The former implies that $\alg{A}\in Ind(\mathsf{L}\oplus\Sigma)$, which together with the latter fact implies that $\Gamma/^{2}\phi\notin \Pi_{2}(Ind(\mathsf{L}\oplus\Sigma))$.
\end{proof}

We then have the following theorem:

\begin{theorem}
There is a dual isomorphism, $Ind$, between the lattice of $\Pi_{2}$-rule systems, and the lattice of inductive rule classes of $\fancyL$-algebras.
\end{theorem}
\begin{proof}
By Proposition \ref{Model Theoretic Part of completeness} and Corollary \ref{Proof Theoretic Part of Completeness}, we have that $Ind$ is a bijection and that $\Pi_{2}$ is its inverse. Also note that $Ind$ and $\Pi_{2}$ are order-reversing maps. Given $\vdash$ and $\vdash^{*}$, two $\Pi_{2}$-rule systems, note that:
\begin{equation*}
    Ind(\vdash \cap \vdash^{*})= Ind(\vdash)\vee Ind(\vdash^{*}).
\end{equation*}
where $Ind(\vdash)\vee Ind(\vdash^{*})$ as usual will mean $\mathbb{R}^{\mathbb{I}}(Ind(\vdash)\cup Ind(\vdash^{*}))$. Indeed, note that if $\alg{A}\in Ind(\vdash)$, or $\alg{A}\in Ind(\vdash^{*})$, then $\alg{A}\in Ind(\vdash\cap\vdash^{*})$, so the right to left inclusion holds. Now to show the other inclusion, we show that:
\begin{equation*}
    \Pi_{2}(Ind(\vdash)\vee Ind(\vdash^{*}))\subseteq {\vdash\cap \vdash^{*}}.
\end{equation*}
Indeed if $\Gamma/^{2}\phi\notin \vdash\cap \vdash^{*}$, then it does not belong to one of them, say $\vdash$. Let $\alg{A}$ be an algebra validating $\vdash$ such that $\alg{A}$ refutes this rule. Then since $\alg{A}\in Ind(\vdash)$, surely $\alg{A}\in Ind(\vdash)\vee Ind(\vdash^{*})$, so $\Gamma/^{2}\phi$ is not valid in such an inductive class, i.e., $\Gamma/^{2}\phi\notin \Pi_{2}(Ind(\vdash)\vee Ind(\vdash^{*}))$.

Additionally we have that:
\begin{equation*}
    Ind(\vdash \vee \vdash^{*})=Ind(\vdash)\cap Ind(\vdash^{*}),
\end{equation*}
Indeed, if $\alg{A}\in Ind(\vdash)$ and $\alg{A}\in Ind(\vdash^{*})$, note that mimicking the derivation relation using the consequence relation in $\alg{H}$, we have that $\alg{A}\in Ind(\vdash\vee \vdash^{*})$. Now assume that $\alg{A}\in Ind(\vdash\vee \vdash^{*})$. Then surely it will validate all rules in $\vdash$ and all rules in $\vdash^{*}$, so $\alg{H}\in Ind(\vdash)\cap Ind(\vdash^{*})$. This shows the result.\end{proof}

In the next section, making use of the results of this section\footnote{We briefly remark that the former arguments and techniques can likewise obtain for us a dual isomorphism between the lattice of $\Pi_{2}$-\textit{multi-conclusion} rule systems and general inductive rule classes algebras (i.e., classes closed under ultraproducts and $\forall$-subalgebras). Whilst we do not pursue this here, we believe such a concept would be fruitful in furthering some connections with the model theory of Boolean and Heyting algebras.}, we will outline some basic notions of admissibility which make sense in the setting of $\Pi_{2}$-rule systems.

\section{Admissibility of Inductive Rules}\label{Admissibility of Inductive Rules}

As noted in the introduction, a lot of the interest in $\Pi_{2}$-rules has come from their usage in axiomatisation of logical systems. Therefore, the question of their admissibility is very natural. One can thus formulate several general questions concerning this:
\begin{itemize}
\item When is an inductive rule admissible?
\item When is an admissible inductive rule derivable?
\item How does this relate to admissibility and structural completeness in the classical case?
\end{itemize}

The above asks whether there are natural analogues of all of these concepts for the case of $\Pi_{2}$-rules. Keeping in mind the usual notion of admissibility from algebraic logic (see e.g. \cite{bergman1988,Rybakov1997-or,Metcalfe2012,Iemhoff2015}), and again fixing a given logical calculus of interest, we can give the following definition:

\begin{definition}
    Let $\Gamma/^{2}\psi$ be a $\Pi_{2}$-rule, and $\mathsf{L}\oplus \Sigma$ some $\Pi_{2}$-rule system. We say that the rule $\Gamma/^{2}\psi$ is \textit{admissible} in $\mathsf{L}\oplus \Sigma$ if for all $\chi$:
    \begin{equation*}
        \vdash_{\mathsf{L}\oplus\Sigma\oplus \Gamma/^{2}\psi}\chi \implies \vdash_{\mathsf{L}\oplus \Sigma}\chi
    \end{equation*}
    We say that the rule is \textit{derivable} if $\Gamma/^{2}\psi\in \mathsf{L}\oplus \Sigma$.
\end{definition}

Note that every derivable rule is admissible. The next Lemma establishes a first semantic criterion to recognise admissible rules, even if a very general and crude one:

\begin{lemma}\label{Lemma connecting syntactic and semantic admissibility}
Let $\mathbf{K}$ be an inductive rule class of $\fancyL$-algebras, and let $\mathsf{L}\oplus \Sigma$ be its dual $\Pi_{2}$-rule system. For each $\Gamma/^{2}\psi$, we have that $\Gamma/^{2}\psi$ is admissible over $\mathsf{L}\oplus \Sigma$ if and only if for each substitution $\sigma$ leaving $B_{c}(\Gamma)$ fixed, $\mathbf{K}\vDash \sigma(\Gamma)$, then $\mathbf{K}\vDash \sigma(\psi)$.
\end{lemma}
\begin{proof}
Assume that $\Gamma/^{2}\psi$ is admissible. Assume that $\sigma$ is some substitution such that $\mathbf{K}\vDash \sigma(\Gamma)$. In other words, if $\lambda_{1},...,\lambda_{n}$ are the terms occuring in $\Gamma$, $\mathbf{K}\vDash \sigma(\lambda_{1}\wedge...\wedge \lambda_{n})$, where the latter is understood as a formula. Hence by completeness, $\vdash_{S}\sigma(\lambda_{1}\wedge...\wedge \lambda_{n})$, and so a fortriori, $\vdash_{S\oplus \Gamma/^{2}\psi}\sigma(\lambda_{1}\wedge...\wedge \lambda_{n})$, where the bound context is left fixed by $\sigma$; using the rule, we have that $\vdash_{S\oplus \Gamma/^{2}\phi}\sigma(\psi)$; so by admissibility, $\vdash_{S}\sigma(\psi)$, which means that $\mathbf{K}\vDash \sigma(\psi)$.

For the other direction, asssume that $\vdash_{\mathsf{L}\oplus \Sigma\oplus\Gamma/^{2}\phi}\psi$. We wish to show that then $\vdash_{\mathsf{L}\oplus \Sigma}\psi$ as well. Using completeness, we have that $\mathbf{K}\vDash \mathsf{L}\oplus \Sigma$ certainly, and we can preserve the relation that $\mathbf{K}$ validates a set of formulas through the derivation. The fact that whenever $\bf{K}\vDash \sigma(\Gamma)$ then $\bf{K}\vDash \sigma(\psi)$ thus works to mimick the applications of $\Gamma/^{2}\psi$, which allows us to conclude, by induction, that $\mathbf{K}\vDash \psi$, i.e., $\vdash_{\mathsf{L}\oplus \Sigma}\psi$.\end{proof}

Because $\Pi_{2}$-rules have a higher complexity, there are alternative notions of admissibility that also make themselves available in this context. Namely, the following makes sense here:

\begin{definition}
Let $\Gamma/^{2}\phi$ be a $\Pi_{2}$-rule and $\mathsf{L}\oplus \Sigma$ some $\Pi_{2}$-rule system. We say that the rule $\Gamma/^{2}\phi$ is \textit{hereditarily admissible} if for all $\Delta/\psi$ a $\Pi_{1}$-rule, we have that
\begin{equation*}
\Delta\vdash_{\mathsf{L}\oplus\Sigma\oplus \Gamma/^{2}\phi}\psi \implies \Delta\vdash_{\mathsf{L}\oplus\Sigma}\psi.
\end{equation*}
\end{definition}

Note that a derivable rule is hereditarily admissible: whenever $\Delta/\psi$ can be derived by using the derived rule, we can carry out the derivation of $\Gamma/^{2}\phi$ to obtain the conclusion without relying on the rule. Hereditary admissibility is moreover, a priori, a stronger property than admissibility: every $\Pi_{1}$-rule which is hereditarily admissible is derivable by definition.

\begin{lemma}\label{Hereditary admissibility and its relation to quasivarieties}
Let $\mathbf{K}$ be an inductive rule class, and $\mathsf{L}\oplus \Sigma$ the dual $\Pi_{2}$-rule system. For each $\Gamma/^{2}\phi$ we have that $\Gamma/^{2}\phi$ is hereditarily admissible if and only if for each $\Delta$ a set of formulas, and $\sigma$ a substitution leaving the bound context fixed, $\mathbf{K}\vDash \sigma[\Delta]\rightarrow \sigma(\psi_{i})$ for each $\psi_{i}\in \Gamma$, then $\mathbf{K}\vDash \sigma[\Delta]\rightarrow \sigma(\phi)$.
\end{lemma}
\begin{proof}
Assume that $\Gamma/^{2}\phi$ is hereditarily admissible. Assume that $\Delta$ is arbitrary and $\sigma$ is some substitution such that $\mathbf{K}\vDash \sigma[\Delta]\rightarrow \sigma(\psi_{i})$ for each $\psi_{i}$. Hence we have that $\sigma[\Delta]\vdash_{\mathsf{L}\oplus \Sigma}\sigma(\psi_{i})$, through a substitution which preserves the bound context of $\Gamma$. Hence surely $\sigma[\Delta]\vdash_{\mathsf{L}\oplus \Sigma\oplus \Gamma/^{2}\phi}\sigma(\psi_{i})$; by using the rule, we then have that $\sigma[\Delta]\vdash_{\mathsf{L}\oplus \Sigma\oplus\Gamma/^{2}\phi}\sigma(\phi)$, and by hereditary admissibility, then $\sigma[\Delta]\vdash_{\mathsf{L}\oplus \Sigma} \sigma(\phi)$, hence $\mathbf{K}\vDash \sigma[\Delta]\rightarrow \sigma(\phi)$.

Now assume the condition. Suppose that $\Delta\vdash_{\mathsf{L}\oplus \Sigma\oplus\Gamma/^{2}\phi}\psi$. Assume that $\psi_{1},...,\psi_{n}$ is a derivation. Then find $\mathcal{A}\in \bf{K}$ which validates the axioms of $\Delta$; then $\mathcal{A}$ can reproduce the derivation given; admissibility is used to witness the use of the rule $\Gamma/^{2}\phi$. Hence we conclude with $\mathbf{K}\vDash \Delta/\psi$, so $\Delta\vdash_{\mathsf{L}\oplus \Sigma}\psi$.
\end{proof}

As promised we can now relate the concepts of admissibility and hereditary admissibility to generation by varieties and quasivarieties:

\begin{lemma}\label{Admissibility by Generation}
Let $\bf{K}$ be an inductive rule class. Suppose that $\mathbb{V}(\bf{K})=\mathbb{V}(\bf{S})$ where $\bf{S}$ is a class of $\fancyL$-algebras validating a rule $\Gamma/^{2}\phi$. Then $\Gamma/^{2}\phi$ is admissible in $\bf{K}$. Similarly if $\bf{K}$ is a quasivariety and $\mathsf{QVar}(\bf{K})=\mathsf{QVar}(\bf{S})$, then $\Gamma/^{2}\phi$ is hereditarily admissible.
\end{lemma}
\begin{proof}
Assume that $\bf{K}\vDash \sigma(\Gamma)$ through a substitution not affecting the bound context of $\Gamma$. Hence $\bf{S}\vDash \sigma(\Gamma)$ as well. Now given the hypothesis, this means that for each $\alg{H}\in \bf{S}$, and each model $v$, $\alg{H},v\vDash \tilde{\forall}\overline{p} \ \Gamma$: indeed, one can simply change the valuation only for propositions occuring in $\overline{p}$, and note that the assumption will then be that $\alg{H},v'\vDash \sigma(\Gamma)$ as well; so by hypothesis, $\alg{H},v\vDash \sigma(\phi)$. So $\bf{S}\vDash \sigma(\phi)$, which immediately implies that $\bf{K}\vDash \sigma(\phi)$.
\end{proof}

We can also provide a form of converse of the previous statement:

\begin{lemma}\label{Admissibility implies generation by elements}
    Let $\Gamma/^{2}\psi$ be a rule which is admissible in an inductive rule class $\bf{K}$. Then $\mathbb{V}(\bf{K})$ is generated as a variety by:
    \begin{equation*}
        \{\alg{H}\in \bf{K} : \alg{H}\vDash \Gamma/^{2}\psi\}
    \end{equation*}
    Similarly if $\Gamma/^{2}\psi$ is hereditarily admissible in an inductive rule class $\bf{K}$ then $\mathsf{QVar}(\bf{K})$ is generated as a quasivariety by the same set.
\end{lemma}
\begin{proof}
Consider:
    \begin{equation*}
        \mathbb{V}(\{\alg{H}\in \bf{K} : \alg{H}\vDash \Gamma/^{2}\psi\})
    \end{equation*}
    We wish to show that this is the same as $\mathbb{V}(\bf{K})$. One inclusion is clear. For the other, assume that $\{\alg{A}\in \bf{K} : \alg{A}\vDash \Gamma/^{2}\psi\}\vDash \chi$ for some formula $\chi$. Then note that:
    \begin{equation*}
        \vdash_{\Pi_{2}(\bf{K})\oplus \Gamma/^{2}\psi}\chi,
    \end{equation*}
    So by admissibility, we have that $\vdash_{\Pi_{2}(\bf{K})}\chi$, i.e., $\bf{K}\vDash \chi$. By Birkhoff's theorem this shows that the two classes are equal. The result for hereditary admissibility and quasivarieties is similar.
\end{proof}

Some logical systems have the desirable property that all admissible rules are derivable -- such systems are often called \textit{structurally complete}. We can also formulate an analogous concept to this one for the setting of $\Pi_{2}$-rule systems:

\begin{definition}
    Let $\mathbf{K}$ be an inductive rule class. We say that $\mathbf{K}$ is \textit{inductively complete} if whenever $\Gamma/^{2}\phi$ is an admissible rule in the corresponding calculus, then it is derivable.
\end{definition}

\begin{lemma}\label{Lemma on admissibility of free algebra}
    Let $\mathbf{K}$ be an inductive rule class, and let $\mathbf{F}_{\omega}(\bf{K})$ be the free algebra on $\omega$-many generators of this class. For each rule $\Gamma/^{2}\psi$, $\Gamma/^{2}\psi$ is admissible over $\bf{K}$ if and only if it is valid in $\mathbf{F}_{\omega}(\bf{K})$.
\end{lemma}
\begin{proof}
First assume that $\Gamma/^{2}\psi$ is admissible over $\bf{K}$. Let $v$ be a valuation over $\bf{F}_{\omega}(\bf{K})$ such that each variant of it in the bound context of $\Gamma$ evaluates to $1$. Then note that we can turn $v$ into a substitution $\sigma$ which leaves the bound context of $\Gamma$ fixed, and does not induce any clash of variables, and such that $\bf{F}_{\omega}(\bf{K})\vDash \sigma(\Gamma)$. Hence $\bf{K}\vDash \sigma(\Gamma)$, so because the rule is admissible, by Lemma \ref{Lemma connecting syntactic and semantic admissibility}, $\bf{K}\vDash \sigma(\psi)$, i.e., $(\bf{F}_{\omega}(\bf{K}),v)\vDash \psi$. This shows the rule is derivable.

For the converse, note that if $\bf{F}_{\omega}(\bf{K})$ makes the rule valid, then since $\mathbb{V}(\bf{K})=\mathbb{V}(\bf{F}_{\omega}(\bf{K}))$, then any derivable rule in $\bf{F}_{\omega}(\bf{K})$ will be admissible in $\bf{K}$ by Lemma \ref{Admissibility by Generation}.
\end{proof}

Using this we easily obtain a characterization of inductive completeness: a class is inductively complete if and only if it is generated as an inductive rule class by the free algebra. Using this, we can obtain a more useful characterization of inductive completeness, which applies when we assume the variety generated by the class to be itself structurally complete:

\begin{lemma}\label{Strong criterion of equivalence between structural completeness and subinductive classes}
Let $\mathbf{K}$ be an inductive rule class, and assume that $\mathbb{V}(\mathbf{K})$ is structurally complete, i.e., $\mathbb{V}(\bf{K})=\mathsf{QVar}(\bf{F}_{\omega}(\bf{K}))$. Then $\bf{K}$ is inductively complete if and only if every proper inductive subclass generates a proper subquasivariety of the quasivariety generated by $\bf{K}$.
\end{lemma}
\begin{proof}
First assume that the inductive rule class is inductively complete, and hence by Lemma \ref{Lemma on admissibility of free algebra}, it is generated by the free algebra. Let $\mathbf{K}'$ be a subinductive class of $\mathbf{K}$, and assume that $\mathbb{S}(\mathbf{K}')=\mathbb{S}(\mathbf{K})$. Then $\mathbf{K}'\subseteq \mathbf{K}\subseteq \mathbb{S}(\mathbf{K}')$, so clearly $\mathbf{K}\subseteq \mathbb{V}(\mathbf{K}')$ so by the Corollary 11.10 of \cite{BurrisSankappanavar}, we have that $\mathbf{F}_{\omega}(\mathbf{K})=\mathbf{F}_{\omega}(\mathbf{K}')$. But since $\mathbf{F}_{\omega}(\mathbf{K}')\in K'$, we know that $\mathbf{F}_{\omega}(\mathbf{K})\in \mathbf{K}'$, which means that $\mathbf{K}\subseteq \mathbb{R}^{\mathbb{I}}(\mathbf{F}_{\omega}(\mathbf{K}))\subseteq \mathbf{K}'$; in other words, $\mathbf{K}'=\mathbf{K}$.

For the other direction, assume that $\mathbf{K}$ is not generated as a subinductive class by $\mathbf{F}_{\omega}(\mathbf{K})$; then $\mathbb{R}^{\mathbb{I}}(\mathbf{F}_{\omega}(\mathbf{K}))$ is a proper inductive subclass. Since $\mathbb{S}(\mathbf{K})$ is structurally complete, $\mathbb{S}(\mathbf{K})=\mathsf{QVar}(\mathbf{F}_{\omega}(\mathbf{K}))$, and surely this will be the same as $\mathbb{S}(\mathbb{R}^{\mathbb{I}}(\mathbf{F}_{\omega}(\mathbf{K}))$. So this is a proper subinductive class which does not generate a proper subquasivariety.
\end{proof}

We will have opportunity to see that these notions can be fruitfully analysed in specific cases. However, in this section, we wish to draw one more relevant connection, namely to the theory of model companions and model completions. The next theorem is analogous to \cite[Theorem 5.4]{Bezhanishvili2022-if}, and gives us a very useful criterion for hereditary admissibility.

\begin{theorem}\label{Embeddability criterion for hereditary admissibility}
Let $\bf{K}$ be an inductive rule class, and $\Gamma/^{2}\psi$ a $\Pi_{2}$-rule. Then $\Gamma/^{2}\psi$ is hereditarily admissible in $\bf{K}$ if and only if every $\alg{A}\in \bf{K}$ embeds into $\alg{B}\in \bf{K}$ such that $\alg{B}\vDash \Gamma/^{2}\psi$.
\end{theorem}
\begin{proof}
First assume that $\Gamma/^{2}\psi$ is hereditarily admissible in $\bf{K}$. Let $\alg{A}$ be some algebra in $\bf{K}$. Look at $T(\alg{A})$, the \textit{open diagram} of $\alg{A}$, i.e.
\begin{equation*}
T(\alg{A})=\{\phi : \phi \text{ is atomic or negated atomic, and }  \alg{A}\vDash \phi\}
\end{equation*}

We will show that $T(\alg{A})$ is consistent with the theory of $\bf{K}$ together with $\Gamma/^{2}\psi$. For this purpose, let $S=\{\phi_{0},...,\phi_{n},\neg \psi_{0},...,\neg \psi_{m}\}$ be finitely many formulas from $T(\alg{A})$, such that all $\phi_{i}$ and $\psi_{j}$ are atomic. Then note that for each $i$, it cannot be that
\begin{equation*}
\Pi_{2}(\bf{K})\cup \{\Gamma/^{2}\psi\}\vDash \forall\overline{x}(\phi_{0}(\overline{x})\wedge...\wedge \phi_{n}(\overline{x}) \rightarrow \chi_{i}(\overline{x}))
\end{equation*}
To see why, note that if this was the case, then we could consider $\Delta=\{\phi_{0},...,\phi_{n}\}$ and the rule $\Delta/\chi_{i}$. The above would imply that $\Pi_{2}(\bf{K})\cup\{\Gamma/^{2}\psi\}\vDash \Delta/\chi_{i}$. Then by completeness, this would imply that $\Delta\vdash_{\Pi_{2}(\bf{K})\oplus \Gamma/^{2}\psi}\chi_{i}$; but this latter fact would imply, by hereditary admissibility, that $\Delta\vdash_{\Pi_{2}(\bf{K})}\chi_{i}$. This contradicts the fact that $\alg{A}\nvDash \Delta/\psi_{i}$ and $\alg{A}\vDash \Pi_{2}(\bf{K})$. Hence for each $i$, we can find $\alg{A}_{i}\in \mathbf{K}$ such that $\alg{A}_{i}\vDash \exists \overline{x} (\phi_{0}(\overline{x})\wedge...\wedge \phi_{n}(\overline{x})\wedge \neg \psi_{i}(\overline{x}))$. Like in Theorem \ref{Collective theorem for special horn classes}, we can then consider $\alg{A}_{0}\times...\times \alg{A}_{n}$, which yields the desired witness to the consistency of $S$.

Hence by compactness we have that we can find $\alg{B}\in \bf{K}$ such that $\alg{B}\vDash \Gamma/^{2}\psi$, and $\alg{B}\vDash T(\alg{A})$; but note that satisfying $T(\alg{A})$ implies that there is an embedding from $\alg{A}$ to $\alg{B}$.

Now assume that the condition holds. We will use Lemma \ref{Hereditary admissibility and its relation to quasivarieties}. Assume that $\bf{K}\vDash \sigma[\Delta]\rightarrow \sigma(\psi_{i})$ for each $\psi_{i}\in \Gamma$. Assume that $\alg{A}\in \bf{K}$, and $\alg{A}\nvDash \sigma[\Delta]\rightarrow \sigma(\phi)$. Then for some model $v$, $(\alg{A},v)\vDash \sigma[\Delta]$, but $(\alg{A},v)\nvDash \sigma(\phi)$. By assumption, there is some $\alg{B}$ which validates the rule in $\bf{K}$ into which $\alg{A}$ embeds. This means that $\alg{B}\nvDash \sigma[\Delta]\rightarrow \sigma(\phi)$ as well. By assumption on $\bf{K}$, we then have that $(\alg{B},v)\vDash \sigma(\psi_{i})$ for each $\psi_{i}$. Since the substitution leaves the bound context invariant, and $\alg{B}$ validates the rule, then $\alg{B},v\vDash \sigma(\phi)$, a contradiction.
\end{proof}

The previous theorem implies that the $\Pi_{2}$-rules holding in a system are exactly those which are valid in the model completion, whenever this exists:

\begin{theorem}\label{Theorem connecting hereditary admissibility and model companions}
Let $\bf{K}$ be an inductive rule class, and suppose that it has a model companion $T^{*}$. Then a $\Pi_{2}$-rule $\Gamma/^{2}\psi$ is hereditarily admissible in $\bf{K}$ if and only if $T^{*}\vDash \chi(\Gamma,\psi)$, the first-order translation of the rule.
\end{theorem}
\begin{proof}
First assume that $T^{*}\vDash \chi(\Gamma,\psi)$. By general model theoretic facts, every algebra in $\bf{K}$ embeds into an existentially closed extension, which will be a model of $T^{*}$. Hence by Theorem \ref{Embeddability criterion for hereditary admissibility}, the rule $\Gamma/^{2}\psi$ is hereditarily admissible.

Conversely, suppose that $\Gamma/^{2}\psi$ is hereditarily admissible in $\bf{K}$. By the above theorem, we have that every algebra in $\bf{K}$ embeds into some algebra satisfying $\Gamma/^{2}\psi$. Let $\alg{A}\vDash T^{*}$. Let $f:\alg{A}\to \alg{B}$ be the extension satisfying $\Gamma/^{2}\psi$. If $\alg{B}\vDash  \exists \overline{y} \phi(\overline{a},\overline{y})$ for $\overline{a}\in \alg{A}$, and $\phi$ an atomic formula, since $\alg{A}$ is existentially closed, then $\alg{A}\vDash \exists\overline{y}\phi(\overline{a},\overline{y})$, as well. Hence the since $\alg{A}$ is existentially closed, being a model of the model companion, it also satisfies $\exists\overline{y}\phi(\overline{a},\overline{y})$. Thus $f$ is a $\forall$-embedding, and thus, by preservation, $\alg{A}\vDash \Gamma/^{2}\psi$.\end{proof}

We also note the following relevant idea:

\begin{proposition}
    Let $\bf{K}$ be an inductive rule class. Then $\Pi_{2}(\bf{K})$ has a model companion if and only if $\mathsf{Log}(\mathsf{QVar}(\bf{K}))$ has a model companion.
\end{proposition}
\begin{proof}
    Assume that $\mathsf{Log}(\mathsf{QVar}(\bf{K}))$ has a model companion, say a theory $T^{*}$. Certainly every model of $\bf{K}$ is a model of $\mathsf{QVar}(\bf{K})$, and so can be embedded in a model of $T^{*}$; conversely, every model of $T^{*}$ can be embedded in a model of the generated quasivariety, which in turn embeds into a model of $\bf{K}$. Hence $T^{*}$ is a model companion of $\Pi_{2}(\bf{K})$.

    Conversely, assume that $\Pi_{2}(\bf{K})$ has a model companion $T^{*}$. Then every model of $T^{*}$ can be embedded in a model of $\mathsf{QVar}(\bf{K})$; additionally, every model of the quasivariety embeds into a model of $\bf{K}$, and through that, into a model of $T^{*}$. So $T^{*}$ is a model companion of the quasivariety. 
\end{proof}

In the next section we will give some examples of how to use these results in a concrete setting. 

\section{Inductive Rule Classes of G\"{o}del algebras}\label{Inductive rule classes of godel algebras}

In these last two sections we provide some examples of how the general theory from before can be fruitfully analysed in specific cases. Throughout fix $\mathsf{IPC}$ as our base calculus. We recall that the \textit{G\"{o}del-Dummett logic} is axiomatised over this system as
\begin{equation*}
    \mathsf{LC}\coloneqq \mathsf{IPC}\oplus (p\rightarrow q)\vee (q\rightarrow p).
\end{equation*}

Such a system was considered by G\"{o}del and Dummett, and extensively studied from thereon (see e.g. \cite{Aguilera2016} for a recent survey containing several open problems in the area). They are algebraized by \textit{G\"{o}del algebras}, with the same axiom defining them. An especially important case of such algebras are the \textit{linear algebras}, which are G\"{o}del algebras which lattice reduct is a linear order; throughout we also refer to these as ``linear Heyting algebras". When these algebras have a second greatest element they are precisely the G\"{o}del algebras which are subdirect irreducibles.

For the sequel, given $n\in \omega$, denote by $[n]$ the unique finite linear G\"{o}del algebra with $n$ many elements; denote by $\lambda_{n}$ the formula which defines the variety generated by this algebra. A fact we will need about $\lambda_{n}$ is that for each $m,k\in \omega$, we have $m\leq k$, if and only if $[m]\sat \lambda_{k}$.

\subsection{Inductive Completeness and Cardinality of Inductive Rule Classes of G\"{o}del Algebras}

Using the results from Section \ref{Admissibility of Inductive Rules} we can give a particularly simple proof of the admissibility of some such rules which have been considered in the context of first-order G\"{o}del logic \cite{Beckmann2008}:

\begin{example}
Recall from Example \ref{Takeuti-Titani Rule} the Takeuti-Titani ``density rule". As mentioned before, such a rule is valid in a linear Heyting algebra if and only if it is dense. So by our above remarks, we obtain:

\begin{corollary}\label{Density rule is admissible}
The density rule is admissible over $\mathsf{LC}$.
\end{corollary}
\begin{proof}
One simply notes that the variety of all G\"{o}del algebras is generated by the dense linear Heyting algebras, since every finite linear Heyting algebra embeds into a dense one. So by Lemma \ref{Admissibility by Generation} we have the result.
\end{proof}
\end{example}

This has the following consequence:

\begin{corollary}
    The class of G\"{o}del algebras is not inductively complete.
\end{corollary}
\begin{proof}
It suffices to consider a linear Heyting algebra which does not validate the Takeuti-Titani; for instance $[3]$.\end{proof}

We can contrast the former with the proof provided in \cite{metcalfemontagnasubstructuralfuzzy}, which involved a delicate syntactic argument. This indicates that the tools developed in this chapter might be useful in carrying out proofs of admissibility  in other more general settings where, like in $\mathsf{LC}$, there is a good grasp of the subdirect irreducibles.

Having a theory to discuss $\Pi_{2}$-rule systems, it is natural to ask what is the structure of such systems in the setting of extensions of $\mathsf{LC}$, which as noted above are quite tame. Indeed, as just mentioned, it is well-known that there are only countably many G\"{o}del logics, all of which are structurally complete \cite{dzikwronskistructuralcompleteness}, and hence, all G\"{o}del quasivarieties are already varieties. Additionally, there are also only countably many first-order G\"{o}del logics \cite{Beckmann2008}. As it turns out, though we do not focus on this setting in the present paper, there are also only countably many multiple conclusion consequence relations\footnote{This result, which is not currently available in the literature, follows from some work communicated to the author by Ian Hodkinson. The key idea lies in showing that the order on co-trees given by surjective p-morphism (i.e., given trees $T,Q$, $T\leq Q$ if there is a surjective p-morphism from $Q$ to $T$) is a well-partial order, which can be done by adapting the classical proof of Kruskal's theorem.}.

In this section we will show that the situation for inductive classes changes dramatically. Given a subset $X\subseteq \omega$, write $\mathbb{R}^{\mathbb{I}}(X)$ for the inductive rule class generated by all $[k]$ for $k\in X$.

\begin{lemma}\label{Separation Lemma for classes of linear inductive classes}
Let $X,Y\subseteq \omega$ be two subsets, and assume that $X\neq Y$. Then $\mathbb{R}^{\mathbb{I}}(X)\neq\mathbb{R}^{\mathbb{I}}(Y)$. Indeed, if $n\notin X$ then $[n]\notin \mathbb{R}^{\mathbb{I}}(X)$.
\end{lemma}
\begin{proof}
Assume that $n\notin X$, but $[n]\in \mathbb{R}^{\mathbb{I}}(X)$. By our Theorem \ref{Collective theorem for special horn classes}, giving the model theoretic completeness, then $[n]\leq_{\forall} \alg{P}$ where $\alg{P}=\prod_{i\in I}\alg{M}_{i}$ where $\alg{M}_{i}$ are ultraproducts of algebras in $X$. Now, by assumption, $[n]$ satisfies $\lambda_{n}$. Hence because this is a $\forall$-subalgebra, $\alg{P}\sat \lambda_{n}$ as well, and so will each of its factors, by preservation of equations in products. Hence by definition of satisfaction of a formula in an ultraproduct, given any $\alg{M}_{i}$ there must be ultrafilter many coordinates where the defining equation is valid; hence, for ultrafilter many coordinates, the algebra at that coordinate is smaller than or equal to $[n]$. 

Now, since $n\notin X$, for each $M_{i}$, there must be a greatest $k<n$ such that $k\in X$ (and this must be non-empty, by the assumption on the ultrafilter), and for ultrafilter many coordinates $m$, $\alg{M}_{i}(m)=[k]$. Hence by Los' theorem, we have that $\alg{M}_{i}\vDash \lambda_{k}$.

Since this argument holds regardless of $\alg{M}_{i}$, we have that $\alg{P}\vDash \lambda_{k}$ as well. But then we have that since $[n]$ is a subalgebra of $P$, then $[n]\sat \lambda_{k}$, from which we obtain that $n\leq k$ -- a contradiction. Hence by reductio, we have that $[n]\notin \mathbb{R}^{\mathbb{I}}(X)$. The statement concerning $X\neq Y$ follows immediately.\end{proof}

\begin{corollary}
There are at exactly $2^{\aleph_{0}}$ inductive rule classes of G\"{o}del algebras.
\end{corollary}
\begin{proof}
The former gives us a lower bound, which is an upper bound since inductive rule classes are bounded by the size of the language.\end{proof}

The former result should sound quite strange indeed if one is used to working with varieties and quasivarieties of Heyting algebras: it implies that with respect to inductive rule classes, the order which makes classical logic the unique maximal extension disappears: the inductive rule class generated by $[2]$ is incomparable with the one generated by $[3]$. Let us look at this situation further. Throughout we write $a/b$ for $(b\rightarrow a)\rightarrow b$. This connective is frequently used in the literature on G\"{o}del algebras, for example by \cite{Baazcompactpropositionallogics}. We note the following about $a/b$, which is easily shown.

\begin{lemma}
For each $a,b\in \alg{H}$ where $\alg{H}$ is a linear Heyting algebra, we have that:
\begin{itemize}
    \item If $a<b$, then $a/b=1$;
    \item If $b\leq a$, then $a/b=b$.
\end{itemize}
\end{lemma}

To give some logical substance to the phenomenon just outlined of the classes being separate, we will give some examples of inductive rules which separate these classes. For instance, consider the following\footnote{We note that strictly speaking this rule is not a $\Pi_{2}$-rule, seeing that it involves $\leq$; but over Heyting algebras this is known to be equivalent, since any such rule can be turned into an equality}:
\begin{equation*}
    \forall p( \forall q(\neg\neg q\leq q\vee p) \rightarrow p=1)
\end{equation*}
Notice that in $[2]$ this rule fails: we have that $0\neq 1$, but if $q$ is arbitrary, then $\neg\neg q=q\leq q\vee p$. So the rule fails in this algebra. On the other hand, it is valid in $\mathbf{3}$: if $p\neq 1$, then let $a$ be the intermediate element, and let $q=a$. Then we have that $\neg\neg q=1$, whilst $q\vee p\leq a$, so the rule holds.

Also note that $\neg\neg q\approx \neg q\rightarrow q$ for linear Heyting algebras, which can be written in our notation as $0/q$.  More generally, consider the following, which are patterned on rules from \cite{Baazcompactpropositionallogics}.
\begin{equation*}
    \rho_{n}\coloneqq \forall p(\forall p_{1},...,p_{n}(0/p_{1}\wedge p_{1}/p_{2}\wedge...\wedge p_{n-1}/p_{n}\leq p_{1}\vee...\vee p_{n}\vee p) \rightarrow p=1)
\end{equation*}
Then we claim that such a rule is valid in $[n+1]$ but not $[n]$. To see why it is not valid in $[n]$, note that given $p\neq 1$, in order for the inequality not to hold, then $p_{i+1}\rightarrow p_{i}$ can never be $1$; hence $p_{i}<p_{i+1}$, which is impossible since this would draw a chain of $n+1$ many elements in $[n]$.

On the other hand, $[n+1]$ validates this: given $p\neq 1$, simply pick precisely the chain of $n$ many elements from $0$ to the element covered by $1$. Then $p_{i+1}\rightarrow p_{i}$ will always be $p_{i}$, so $(p_{i+1}\rightarrow p_{i})\rightarrow p_{i+1}$ will be $1$, this holding for all terms on the left; whilst on the right no element is $1$. More generally, any chain of size greater than or equal to $n$ will validate this rule.

\subsection{$\forall$-subalgebras of G\"{o}del algebras}

Whilst in general the notion of a $\forall$-subalgebra appears to be difficult to outline, in the context of \textit{linear} Heyting algebras we can provide a more concrete description of such embeddings. Recall that given any partial order $(P,\leq)$ we write $a\prec b$ if and only if $a<b$, and whenever $a<c$, then $b\leq c$.

\begin{definition}
    Let $f:\alg{H}\to \alg{H}'$ be an embedding of linear Heyting algebras. We say that $f$ is \textit{cover-preserving} if whenever $a\prec b$ then $f(a)\prec f(b)$.
\end{definition}

\begin{theorem}
Let $f:\alg{H}\to \alg{H}'$ be an embedding of linear Heyting algebras. Then $f$ is a $\forall$-embedding if and only if it is cover-preserving.
\end{theorem}
\begin{proof}
First we note that every $\forall$-embedding is cover-preserving. To see this, consider the formula:
\begin{equation*}
    \forall c (a/c \wedge c/b \leq c\vee a)
\end{equation*}
Assume that $a\prec b$ but there is some $c$ such that $f(a)<c$ and $c<f(b)$. Hence note that $f(a)/c=1$ and $c/f(b)=1$, which is not below $f(a)\vee c$; so $f$ could not be a $\forall$-embedding. On the other hand, if $c$ is arbitrary in $\alg{H}$; then $c\leq a$ or $b\leq c$, and one can verify that the above equation holds. Hence since $f$ is a $\forall$-embedding, we obtain that $f(a)\prec f(b)$.

Now assume that $f$ is a cover-preserving embedding. Let $\overline{a}$ be some collection of elements in $\alg{H}$; we will show that $(\alg{H}',f\overline{a})\in \mathbb{SP}_{U}(\alg{H},\overline{a})$. This amounts to showing, by a classical result, that every finite subgraph can be embedded into $(\alg{H},\overline{a})$. So let $\mathbb{X}$ be such a local subgraph of $(\alg{H},f\overline{a})$, and assume without loss of generality that it contains all constants.

Now define an embedding from $\mathbb{X}$ to $(\alg{H},\overline{a})$ as follows: first, send $f(a)$ to $a$; now for each other $x\in X$, notice that we can identify this element, in $\alg{H}'$, by its position relative to all elements previously identified. We order the constants as $0,a_{0},...,a_{n},1$, and proceed to define the mapping within the blocks.

So first, assume that only the constants have been defined. We work between $0$ and $a_{0}$, but the same argument will work wherever. Enumerate the elements ocurring in $0<x<f(a_{0})$, in order, as $x_{0},...,x_{n}$. Notice that if the interval $[0,a_{0}]$ is finite with cardinality $n$, then so is the interval $[0,f(a_{0})]$: formally one can prove this by induction, by noting that each element in $[0,a_{0}]$ will be a cover of the previous element, and using the fact that the map is cover-preserving. Hence if the interval is finite, we have an onto isomorphism, and can map $x_{0},...,x_{n}$ in the same order. If the interval is infinite, then we can always find new images for the elements.

Proceeding in this way, we have that we can define the map $g$ from $\mathbb{X}$ to $\alg{H}'$. Such a map is a local subgraph embedding, since the order is preserved (and hence the meet and join properties are respected, and so is the implication). Hence $(\alg{H}',f\overline{a})\in \mathbb{SP}_{U}(\alg{H},\overline{a})$ and so:
\begin{equation*}
    (\alg{H}',f\overline{a})\in \mathbb{V}(\alg{H},\overline{a})
\end{equation*}
Hence by Proposition \ref{Criterion for being a universal subalgebra}, the map $f$ is a $\forall$-embedding, as desired.\end{proof}

We will have occasion to make use of the above ideas further in the next section. Before that, we make some final notes about $\forall$-embeddings. Readers familiar with Esakia duality will have wondered what the dual of a $\forall$-embedding ought to be. Whilst we are far from a definitive answer, it seems to us that such a condition should be related to the notion of a local homeomorphism as discussed in \cite{kuznetsovlocalhomeomorphisms}.

\section{Inductively Complete Inductive Classes and Admissibility of $\Pi_{2}$-rules}\label{Inductively Complete Inductive Classes and Admissibility of Pi-2 rules}

In this final section we provide two theorems which give us some positive results contrasting to the complex picture painted by the previous section: despite the fact that there are continuum many inductive classes of G\"{o}del algebras, there are only countably many inductively complete ones, which are in 1-1 correspondence with the varieties and quasivarieties. These are precisely the minimal inductive classes. We also prove that the admissibility problem of $\Pi_{2}$-rules for $\mathsf{LC}$ is decidable, following from our work from Section \ref{Admissibility of Inductive Rules}, given that $\mathsf{LC}$ has a model completion.

Both of our results are consequences of the following technical structure theorem. It was pointed out to us that such a theorem might be folklore in the literature on G\"{o}del algebras; we opted to include a full proof of it, for the benefit of the reader. Moreover, we note that $\mathbb{Q}$ below denotes the chain of the rationals with a minimal and maximal point added to it.

\begin{theorem}\label{Every algebra contains some form of chain attached to it}
    Let $\alg{H}$ be a G\"{o}del algebra. Then either $[n]\leq_{\forall}\alg{H}$ or $\mathbb{Q}$ is a $\forall$-subalgebra of an ultrapower of $\alg{H}$.
\end{theorem}
\begin{proof}
    Because $\mathsf{LC}$ is locally finite,  given $\alg{H}$ we have two options:
    \begin{enumerate}
        \item There is a uniform bound $n$ on the size of chain-subalgebras of $\alg{H}$. Then first we claim that $[n]\leq_{\forall}\alg{H}$. It is clear that it is a subalgebra, say through an embedding $f$. Now pick all constants $\overline{n}$, and we need to show that:
        \begin{equation*}
            (\alg{H},f\overline{n})\in \mathbb{V}([n],\overline{n}).
        \end{equation*}
        To see this in turn, pick some finite subalgebra of $\alg{H}$ containing $f\overline{n}$, say $(\alg{L},f\overline{n})$. By universal algebra, we have that
        \begin{equation*}
            (\alg{L},f\overline{n})\leq \prod_{i=1}^{n}\alg{L}_{i}
        \end{equation*}
        through a subdirect inclusion, in the variety $[n]$, with no added constants, where $\alg{L}_{i}$ gets its constants interpreted by forcing their interpretation. Now first notice that $(\alg{L}_{i},f\overline{n}(i))$ can only identify constants with $1$: indeed if $p(f(n))\neq 1$, and $f(n)\neq f(m)$ then without loss of generality, since they are linearly ordered, $f(n)<f(m)$, so $f(m)\rightarrow f(n)=f(n)$. If $p:\alg{L}\to \alg{L}_{i}$ is the homomorphism, then $p(f(m)\rightarrow f(n))=p(f(n))$; but if $p(f(m))=p(f(n))$, then $p(f(m)\rightarrow f(n))=1$, a contradiction.

        Hence the constants will appear in $\alg{L}_{i}$ in the correct order, possibly with some of them identified with $1$. Hence we can take a homomorphic image of $[n]$ which obtains an isomorphic copy to $(\alg{L}_{i},pf\overline{n})$. Hence we obtain that $(\alg{L},\overline{f}\overline{n})$ will be a subdirect product of such homomorphic images, and so it will belong to $\mathbb{V}([n],\overline{n})$. Hence we have that each finite subgraph of $(\mathcal{H},f\overline{n})$ belongs to $\mathbb{V}([n],\overline{n})$, which shows, because the latter is then a universal class, that $(\mathcal{H},f\overline{n})$ belongs there as well.
        \item There is no bound on the chains. Then note that $\mathbb{Q}$ will embed in an ultrapower of $\alg{H}$ by usual model theoretic reasoning, using compactness. Hence $\mathbb{Q}\leq \mathcal{D}$ for some such ultrapower. Then we show that for each collection of constants $f\overline{a}\in \mathbb{Q}$,
        \begin{equation*}
            (\mathcal{D},f\overline{a})\in \mathbb{V}(\mathbb{Q},\overline{a}).
        \end{equation*}
        The argument proceeds as above, by taking a finite subgraph, $(\alg{L},f\overline{a})$, which again can be decomposed into linear factors. The additional argument needed is that $\overline{a}$ forms a linear subalgebra of $\mathbb{Q}$, from which we then proceed to construct a homomorphism; the argument concludes then in the same way.
    \end{enumerate}

\end{proof}

\subsection{Inductively Complete $\mathsf{LC}$-Inductive Rule Classes}

We have shown that $\mathsf{LC}$ is not inductively complete, so a natural question arises whether anything can be said about those extensions of it which are inductively complete. For general logical systems this appears to be a difficult question, but making use of the structural completeness of $\mathsf{LC}$ and the previous structure theorem, we can provide a complete characterisation. This goes through the following lemma:

\begin{lemma}
    The inductive rule class $\mathbb{R}^{\mathbb{I}}(\mathbb{Q})$, and the inductive rule classes $\mathbb{R}^{\mathbb{I}}([n])$ are all minimal; they are precisely the minimal inductive classes.
\end{lemma}
\begin{proof}
Assume that $\bf{K}\subseteq \mathbb{R}^{\mathbb{I}}(\mathbb{Q})$ and this is non-empty. Assume that they are different; then we may assume that $\mathbb{Q}\notin \bf{K}$. Let $\alg{H}$ be some element in $\bf{K}$, which exists by non-emptiness. By the previous theorem, if $\mathbb{Q}$ did not embed in an ultrapower of $\alg{H}$, then $[n]\leq_{\forall}\alg{H}$ for some $n$; but this is a contradiction, since then the equation characterising $[n]$ would be valid on a product of dense linear algebras, a contradiction. So $\mathbb{Q}$ embeds in an ultrapower of $\alg{H}$; but such an ultrapower will also be $\bf{K}$, and hence, so will $\mathbb{Q}$.

Similarly, assume that $\bf{K}\subseteq \mathbb{R}^{\mathbb{I}}([n])$. Again suppose that $[n]\notin \bf{K}$. Let $\alg{H}$ be some algebra in $\bf{K}$. Note that then $\alg{H}$ must have a uniform bound on its chains (otherwise $\mathbb{Q}$ would be in $[n]$, which is absurd). Moreover this uniform bound must be $n$, otherwise by the arguments above, $[m]\leq_{\forall}\alg{H}$ which would belong to $\mathbb{R}^{\mathbb{I}}([m])$, in violation of Lemma \ref{Separation Lemma for classes of linear inductive classes}. Hence by the above result $[n]\leq_{\forall}\alg{H}$, showing that $[n]\in \bf{K}$.

Now assume that $\bf{K}$ is a minimal class (hence non-empty). Let $\alg{H}$ be some algebra in $\bf{K}$. Then by the above theorem, either $[n]\in \bf{K}$ or $\mathbb{Q}\in \bf{K}$. Since $\bf{K}$ is minimal, then it will coincide with one of the former, as desired.
\end{proof}

From this we get a characterization of the inductively complete inductive rule classes:

\begin{theorem}\label{Characterization of the inductively complete rule classes}
    The inductively complete inductive rule classes in $\mathsf{LC}$ are precisely the minimal ones.
\end{theorem}
\begin{proof}
    Throughout we use Lemma \ref{Strong criterion of equivalence between structural completeness and subinductive classes}. Certainly all minimal inductive rule classes are inductively complete, since they do not have any proper subinductive classes. Now let $\bf{K}$ be some inductive rule class which is not minimal. First assume that the algebras in $\bf{K}$ all have a uniform bound of $n$. Then note that $[n]$ is a proper subinductive class (because $\bf{K}$ is not minimal), but the quasivariety generated by it generates the quasivariety generated by $\bf{K}$. Now assume that there is no such bound. Then $\mathbb{Q}\in \bf{K}$, so again we get the same conclusion by looking at the minimal inductive class $\mathbb{R}^{\mathbb{I}}(\mathbb{Q})$. Hence $\bf{K}$ cannot be inductively complete.
\end{proof}

To construct an axiomatisation of $\mathbb{R}^{\mathbb{I}}([n])$, one may extract such an axiomatisation from an axiomatisation of the model companion of $\mathbb{V}([n])$. In the special case of $\mathbb{V}(\mathbb{Q})$, this can be effectively studied using the results from \cite{darnierejunkermodelcompletionofcoheyting}, characterising the model completion of $\mathsf{LC}$, and implies in particular that the $\Pi_{2}$-rule system associated to $\mathbb{R}^{\mathbb{I}}(\mathbb{Q})$ is recursively axiomatisable, and decidable. In both cases, this problem appears to be difficult and intimately related to the study of the model companions and completions of such varieties, and therefore it is left outside of the scope of the present article. 

\subsection{Decidability of Admissibility of $\Pi_{2}$-Rules of $\mathsf{LC}$}

In \cite{Bezhanishvili2022-if} the authors studied several procedures to show the decidability of the problem of admissibility of $\Pi_{2}$-rules over several logical systems. These concerned rules as outlined in Example \ref{Pi2rulesNIckandSilvioversion}, which crucially leaves out settings such as the present one. As it turns out, some of the results presented therein can be generalized. Key to our work will be the connection between hereditary admissibility and admissibility spelled out in Theorem \ref{Theorem connecting hereditary admissibility and model companions}. We start with a decidability result that is somewhat weaker:

\begin{theorem}\label{Effective recognizability of hereditary admissibility}
    Hereditary Admissibility of $\Pi_{2}$-rules over $\mathsf{LC}$ is effectively recognizable.
\end{theorem}
\begin{proof}
By Theorem \ref{Theorem connecting hereditary admissibility and model companions}, this amounts to showing that in the model completion $\mathsf{LC}^{*}$ we can decide whether $\mathsf{LC}^{*}\vDash \Gamma/^{2}\phi$. We know that the theory $\mathsf{LC}^{*}$ will be complete\footnote{Here is a short argument: By exercise 2 of Section 8.3 of Hodges, a first order theory which is model complete and has the joint embedding property is complete. Now $\mathsf{LC}$ certainly has the joint embedding property, and this implies that its model companion does as well.}, which implies that it is decidable. But additionally, quantifier elimination is effectively recognizable by the same arguments as used in \cite[Corollary 5.7, Lemma 5.8]{Bezhanishvili2022-if}. Hence, since $\mathsf{LC}$ is decidable (and therefore we can decide the validity of quantifier-free formulas as computed by our theory), this implies the whole problem is effectively recognizable.
\end{proof}

One immediately obtains that in fact hereditary admissibility of $\Pi_{2}$-rules of a system with a model completion will be effectively recognizable. But one might naturally wonder about admissible rules which are not hereditarily admissible. As we will show, in the case of $\mathsf{LC}$, the former already gives us everything we need, since all admissible rules are hereditarily admissible. To see this, we begin by stating an easy fact.

\begin{proposition}
    A rule $\Gamma/^{2}\phi$ is admissible/hereditarily admissible over $\mathsf{LC}$ if and only if it is admissible/hereditarily admissible over the class of linear Heyting algebras.
\end{proposition}
\begin{proof}
    Obvious, since $\mathsf{LC}$ is generated as an inductive rule class by the linear algebras.
\end{proof}

\begin{proposition}\label{Admissibile but not hereditarily admissible implies only finitely many finite chains}
    Let $\Gamma/^{2}\phi$ be a rule which is not hereditarily admissible over $\mathsf{LC}$. Then
    $\Gamma/^{2}\phi$ is valid in only finitely many chains and no infinite chain.
\end{proposition}
\begin{proof}
If $\Gamma/^{2}\phi$ is not hereditarily admissible, by Theorem \ref{Embeddability criterion for hereditary admissibility}, there is an algebra $\alg{H}$ which is not embeddable in an algebra satisfying this rule. Now if each finite chain-factor of $\alg{H}$ was embeddable in such an algebra, then any product would be embeddable in a product of the targets, and so $\alg{H}$ (which is a subdirect product of its finite chain subalgebras) would also embed there. So not all chains of $\alg{H}$ can embed into such an algebra, say, $[n]$ does not embed into any such algebra. But if $n\leq m$ and $[m]$ can embed there, then so can $[n]$; moreover no infinite chain can embed there, otherwise all finite chains would embed into it as well.
So there must be only finitely many chains which embed into algebras which validate the rule; since the chains embed into themselves, this implies that there are only finitely many chains which validate the rule.
\end{proof}

Now we combine Lemma \ref{Admissibility implies generation by elements}, with the facts proven at the end of the last section:

\begin{proposition}
    Let $\Gamma/^{2}\phi$ be a rule which is admissible over $\mathsf{LC}$. Then $\mathsf{LC}$ is generated as a variety by the chains validating $\Gamma/^{2}\phi$.
\end{proposition}
\begin{proof}
    By Lemma \ref{Admissibility implies generation by elements} we know that $\mathsf{LC}$ is generated by
    \begin{equation*}
        \{\alg{H} : \alg{H}\vDash \Gamma/^{2}\phi\}.
    \end{equation*}
    Now for each such $\alg{H}$ either it has a bound on chains or it does not. If it does, then if $n$ is such a bound, $[n]\leq_{\forall}\alg{H}$, and so $[n]\vDash \Gamma/^{2}\phi$, and certainly $\alg{H}$ can be recovered from $[n]$ through operations generating a variety; if not $\mathbb{Q}\leq_{\forall}\alg{H}$ and $\mathbb{Q}$ certainly also generates $\alg{H}$ in this manner. So this proves the result.
\end{proof}

This has the following strong consequence, where $\mathbf{Lin}$ is the class of linear Heyting algebras:

\begin{proposition}\label{Admissibility implies cofinally many chains generate}
    Let $\Gamma/^{2}\phi$ be a rule which is admissible over $\mathsf{LC}$. Then the class:
    \begin{equation*}
        \mathbf{Lin}\cap \{\alg{H} : \alg{H}\vDash \Gamma/^{2}\phi\}
    \end{equation*}
    must either include $\mathbb{Q}$ or cofinally many finite chains.
\end{proposition}
\begin{proof}
    Indeed, if not, then there would only be finitely many finite chains in such a set. But then the variety generated by this finite set of chains would not be $\mathsf{LC}$.
\end{proof}

\begin{corollary}\label{Admissibility if and only if hereditary admissibility}
    Over $\mathsf{LC}$, a rule is admissible if and only if it is hereditarily admissible.
\end{corollary}
\begin{proof}
    Assume that $\Gamma/^{2}\phi$ is admissible but not hereditarily admissible. By Proposition \ref{Admissibile but not hereditarily admissible implies only finitely many finite chains} then the class of finite chains validating this rule must be finite, and $\mathbb{Q}$ cannot validate it. But then this contradicts Proposition \ref{Admissibility implies cofinally many chains generate}.
\end{proof}

\begin{corollary}\label{Admissibility of Pi2-rules is effectively recognizable}
    Admissibility of $\Pi_{2}$-rules over $\mathsf{LC}$ is effectively recognizable.
\end{corollary}
\begin{proof}
    This follows from Corollary \ref{Admissibility if and only if hereditary admissibility} together with Theorem \ref{Effective recognizability of hereditary admissibility}.
\end{proof}

Note that the previous proof involved a great deal of universal algebraic concepts. Namely:
\begin{enumerate}
    \item Interpolation and Uniform interpolation (in the form of model completions);
    \item Local finiteness of algebras;
    \item A Theorem, like Theorem \ref{Every algebra contains some form of chain attached to it} which entails that generating classes can always be taken from the `tame ones' (which in these case are the linear algebras).
\end{enumerate}

It is thus natural to wonder whether these concepts are enough to ground the arguments at hand. Several questions could be asked in this direction, but we restrict ourselves to one here. We note that in the above proof the structural completeness of $\mathsf{LC}$ is implicitly used. Additionally, it makes conceptual sense that in the presence of structural completeness the problem of hereditary admissibility should resemble admissibility -- after all, only the $\Pi_{1}$-fragment should remain. However, it is not obvious how one would prove such a general result. Nevertheless, the importance of such a concept merits the following conjectures, in increasing degree of generality:

\begin{conjecture}
    Let $L$ be an intermediate logic which is hereditarily structurally complete. Then hereditary admissibility is equivalent to admissibility over $L$.
\end{conjecture}

\begin{conjecture}
    Let $L$ be an intermediate logic which is structurally complete. Then hereditary admissibility is equivalent to admissibility over $L$.
\end{conjecture}

\section{Conclusion}\label{Section: Conclusion}

In this article I have proposed an analytic framework to study $\Pi_{2}$-rules using tools stemming from model theory and universal algebra. In it I introduce the notions of a $\Pi_{2}$-rule system and an inductive rule class, and establish an algebraic completeness connecting the two concepts. As an example use-case, I study G\"{o}del algebras, showing that the present tools allow for easier admissibility proofs for interesting rules like the Takeuti-Titani rule. Additionally, I present a preliminary study of inductive rule classes over $\mathsf{LC}$, emphasising two aspects: the structure of the lattice of inductive rule classes; and admissibility of rules over such classes. Our key contribution in this respect lies in providing a full characterisation of the inductively complete inductive rule classes, as well as showing the problem of admissibility of inductive rule classes to be decidable over $\mathsf{LC}$.

In this respect there are several follow-up problems that stem immediately from the current work. With respect to $\mathsf{LC}$, these concern, for instance, the axiomatisation of the inductively complete rule classes. Given the results presented in this paper, this can be carried out in parallel and in support of a study of the model companions of the finitely generated varieties of G\"{o}del algebras. Similarly, it would be interesting to understand how the structure of such inductive rule classes through a duality-theoretic perspective.

From a theoretical point of view, it would be interesting to understand the connections between the concepts of hereditary admissibility, admissibility, and the theory of model completions, and the role of structural completeness in establishing these connections. Moreover, it would be interesting to understand how such concepts in turn relate to the theory of unification, as pointed out in \cite{Bezhanishvili2022-if}, and whether $\Pi_{2}$-rules can provide logical meaning to natural computational problems of this sort; this has been explored by the author and Ghilardi in \cite{almeidaghilardiunificationsimplevariable}, but much remains to be done. Finally, as a more long term ambition, it would be interesting to provide a systematic definability theory which explains which classes of frames can be captured using $\Pi_{2}$-rules, and providing conditions for such definitions to hold.

\section{Acknowledgements}

I thank the comments  of the two anonymous reviewers, which substantially improved the presentation of the results. I am thankful to
Silvio Ghilardi and Nick Bezhanishvili for suggestions and discussions about the results of this paper.

\printbibliography

\end{document}